\documentclass[11pt,a4paper]{article}
\usepackage{pdfsync} % FK -- added

\usepackage{amsfonts,amsmath,amsthm}
\newtheorem{theorem}{Theorem}
\numberwithin{theorem}{section}
\newtheorem{lemma}[theorem]{Lemma}
\numberwithin{equation}{section}
\usepackage{amssymb}
\usepackage{hyperref}
\usepackage{graphicx}
\usepackage{wrapfig}
\usepackage{subfig}
\usepackage{epsfig}
\usepackage{float,mathrsfs}
\usepackage{bbm}
\usepackage{algpseudocode}
\usepackage[font=scriptsize]{caption}
\usepackage{chngcntr} % Comment this out if you're using the MacTeX distribution! -VK
\counterwithout{figure}{section}
\usepackage{color}
\definecolor{refkey}{rgb}{0.9451,0.2706,0.4941}\definecolor{labelkey}{rgb}{0.9451,0.2706,0.4941}

\definecolor{darkred}{RGB}{139,0,0}
\definecolor{darkgreen}{RGB}{0,100,0}
\definecolor{darkmagenta}{RGB}{139,0,139}

\newcommand{\setu}{{\mathrm{\mathfrak{u}}}}

\newcommand{\bsx}{{\boldsymbol{x}}}
\newcommand{\bsy}{{\boldsymbol{y}}}

\newcommand{\bsz}{{\boldsymbol{z}}}

\newcommand{\bsgamma}{{\boldsymbol{\gamma}}}
\newcommand{\rd}{{\mathrm{d}}}
\newcommand{\bbZ}{{\mathbb{Z}}}
\newcommand{\bbN}{{\mathbb{N}}}

\newcommand{\satop}[2]{\stackrel{\scriptstyle{#1}}{\scriptstyle{#2}}}

\newcommand{\bsnu}{{\boldsymbol{\nu}}}
\newcommand{\bsOmega}{{\boldsymbol{\Omega}}}

\newcommand{\bszero}{{\boldsymbol{0}}}
\newcommand{\bsone}{{\boldsymbol{1}}}

\newcommand{\bsp}{{\boldsymbol{p}}}

\newcommand{\calO}{\mathcal{O}}

\newcommand{\mask}[1]{{}}

\title{Uncertainty quantification using periodic random variables}

\author{V. Kaarnioja\footnotemark[2]
\and F.~Y. Kuo\footnotemark[2]
\and I.~H. Sloan\footnotemark[2]}

\renewcommand{\thefootnote}{\fnsymbol{footnote}}

\begin{document}
\maketitle

\footnotetext[2]{School of Mathematics and Statistics, University of New South Wales, Sydney NSW 2052, Australia ({\tt v.kaarnioja@unsw.edu.au}, {\tt f.kuo@unsw.edu.au}, {\tt i.sloan@unsw.edu.au})}

\begin{abstract}
Many studies in uncertainty quantification have been carried out under the assumption of an input random field in which a countable number of independent random variables are each uniformly distributed on an interval, with these random variables entering \emph{linearly} in the input random field (the so-called affine model). In this paper we consider an alternative model of the random field, in which the random variables have the same uniform distribution on an interval, but the random variables enter the input field as {\it periodic} functions. The field is constructed in such a way as to have the same mean and covariance function as the affine random field. Higher moments differ from the affine case, but in general the periodic model seems no less desirable. The new model of the random field is used to compute expected values of a quantity of interest arising from an elliptic PDE with random coefficients. The periodicity is shown to yield a higher order cubature convergence rate of $\mathcal{O}(n^{-1/p})$ independently of the dimension when used in conjunction with rank-1 lattice cubature rules constructed using suitably chosen {\em smoothness-driven product and order dependent weights}, where $n$ is the number of lattice points and $p$ is the summability exponent of the fluctuations in the series expansion of the random coefficient. We present numerical examples that assess the performance of our method.
\end{abstract}

% Some footnote styles
%\renewcommand{\thefootnote}{\fnsymbol{footnote}}
%\renewcommand{\thefootnote}{\dag}
\renewcommand{\thefootnote}{\arabic{footnote}}

\section{Introduction}

This paper is concerned with the development and use of specially designed
 random fields on a physical domain $D\subseteq\mathbb{R}^d$,
where $d=1$, $2$, or $3$. For simplicity we assume that the boundary $\partial
D$ is Lipschitz.

Many studies in uncertainty quantification are modeled by PDEs over the domain $D$, in which one or more of the
coefficients is a random field over $D$.  In particular,  many recent
papers (including \cite{cds10,dicklegiaschwab,spodpaper14,ghs18,kss12,kssmultilevel,schwab13}) have used an ``affine'' model of the random
field, taking the form
\begin{equation}\label{eq:A_affine}
A(\bsx, \omega) = \overline{a}(\bsx)+\sum_{j\geq 1} Y_j(\omega)       \,\psi_j(\bsx),
\quad \bsx \in D,~\omega \in \Omega,
\end{equation}
where $(\Omega, A, \mathbb{P})$ is a probability space, 
%\begin{equation*}
$Y_1, Y_2,\ldots$ are independently and identically distributed random variables uniformly distributed on 
 $[-\tfrac{1}{2}, \tfrac{1}{2}],$
%\end{equation*}
and $(\psi_j)_{j\geq 1}$ are  real-valued $L_\infty$ functions on $D$ satisfying
\begin{equation}\label{eq:summable} \sum_{j\geq 1}\|\psi_j\|_{L_\infty}<\infty
\end{equation}
and otherwise for the moment not specified. With these definitions the sum
in \eqref{eq:A_affine} converges uniformly on $D$ for all values of the
$Y_j$, and the random field is pointwise well defined.
The expected value at
$\bsx\in D$ is
\begin{align*}
\mathbb{E}[A(\bsx,\cdot)] =\overline{a}(\bsx)+\sum_{j\geq 1} \mathbb{E}[Y_j]\, \psi_j(\bsx)=\overline{a}(\bsx),
\end{align*}
since
$
\mathbb{E}[Y_j] = \int_{-1/2}^{1/2} y\, \rd y= 0$ for $j\ge 1.
$

The fact that the random variable in \eqref{eq:A_affine} occurs linearly
seems to be a result of history, rather than something imposed by modeling
assumptions.  Suppose instead that we replace $Y_j$ in \eqref{eq:A_affine}
by $\Theta(Y_j)$, where $\Theta\!:[-\tfrac{1}{2}, \tfrac{1}{2}]\to
\mathbb{R}$ is a continuous function with the  properties
\begin{equation}\label{eq:alpha2}
\int_{-1/2}^{1/2}\Theta(y)\,{\rm d}y=0\quad\text{and}\quad\int_{-1/2}^{1/2}\Theta^2(y)\,\rd y= \frac{1}{12},
\end{equation}
both of which are satisfied by the special affine choice $\Theta(y)=y$. Thus we replace \eqref{eq:A_affine} by
\begin{equation}\label{eq:A_affine_alpha}
A(\bsx, \omega) = \overline{a}(\bsx)+\sum_{j\geq 1} \Theta(Y_j(\omega))\,\psi_j(\bsx),
\quad \bsx \in D,~\omega \in \Omega.
\end{equation}
It has exactly the same mean as \eqref{eq:A_affine} because $\mathbb{E}[\Theta(Y_j)]=0$ and also has the same covariance
\begin{align}
\mathrm{cov}(A)(\bsx,\bsx')
\,:=\,&\mathbb{E}[(A(\bsx,\cdot)-\overline{a}(\bsx))(A(\bsx',\cdot)-\overline{a}(\bsx'))]\notag\\
=\,& \sum_{j\geq 1} \sum_{j'\geq 1} \mathbb{E}[\Theta(Y_j)\,\Theta(Y_{j'})]\, \psi_j(\bsx)\,\psi_{j'}(\bsx')\notag\\
=\,&\frac{1}{12}\sum_{j\geq 1}   \psi_j(\bsx)\,\psi_j(\bsx'),\label{eq:mercerlike}
\end{align}
because $\mathbb{E}[\Theta(Y_j)\,\Theta(Y_{j'})]$ vanishes for $j\ne j'$ by
the independence of $Y_j$ and $Y_{j'}$ (and hence of $\Theta(Y_j)$ and
$\Theta(Y_{j'})$), and $\mathbb{E}[\Theta^2(Y_j)]=\frac{1}{12}=\mathbb{E}[Y_j^2]$ from~\eqref{eq:alpha2}. In particular, the variance, obtained by setting $\bsx=\bsx'$, is
\[
\mathrm{var}(A)(\bsx)=\frac{1}{12}\sum_{j\geq 1}   \psi_j^2(\bsx),
\]
independently of the choice of the function $\Theta$.

In this paper we explore a periodic choice of $\Theta$ satisfying~\eqref{eq:alpha2}, namely,
\begin{align*}
\Theta(y) = \frac{1}{\sqrt{6}}\sin(2\pi y), \quad  y \in [-\tfrac12,\tfrac12];
\end{align*}
thus our model of the random field, instead of \eqref{eq:A_affine}, becomes
\begin{equation} \label{eq:A_affine_sin}
A(\bsx, \omega) = \overline{a}(\bsx)+\frac{1}{\sqrt{6}}\sum_{j\geq 1} \sin(2\pi Y_j(\omega))\,\psi_j(\bsx),
\quad \bsx \in D,~\omega \in \Omega.
\end{equation}

\begin{figure}[!t]
\centering
\subfloat[An ensemble of realizations drawn from the affine random field $A_{\rm aff}(x,\omega)=\overline{a}(x)+\sum_{j=1}^{100}Y_j(\omega)\,\psi_j(x)$, $x\in[0,1{]}$, $Y_j\sim U([-\frac12,\frac12{]})$.]
{{\includegraphics[height=.36 \textwidth]{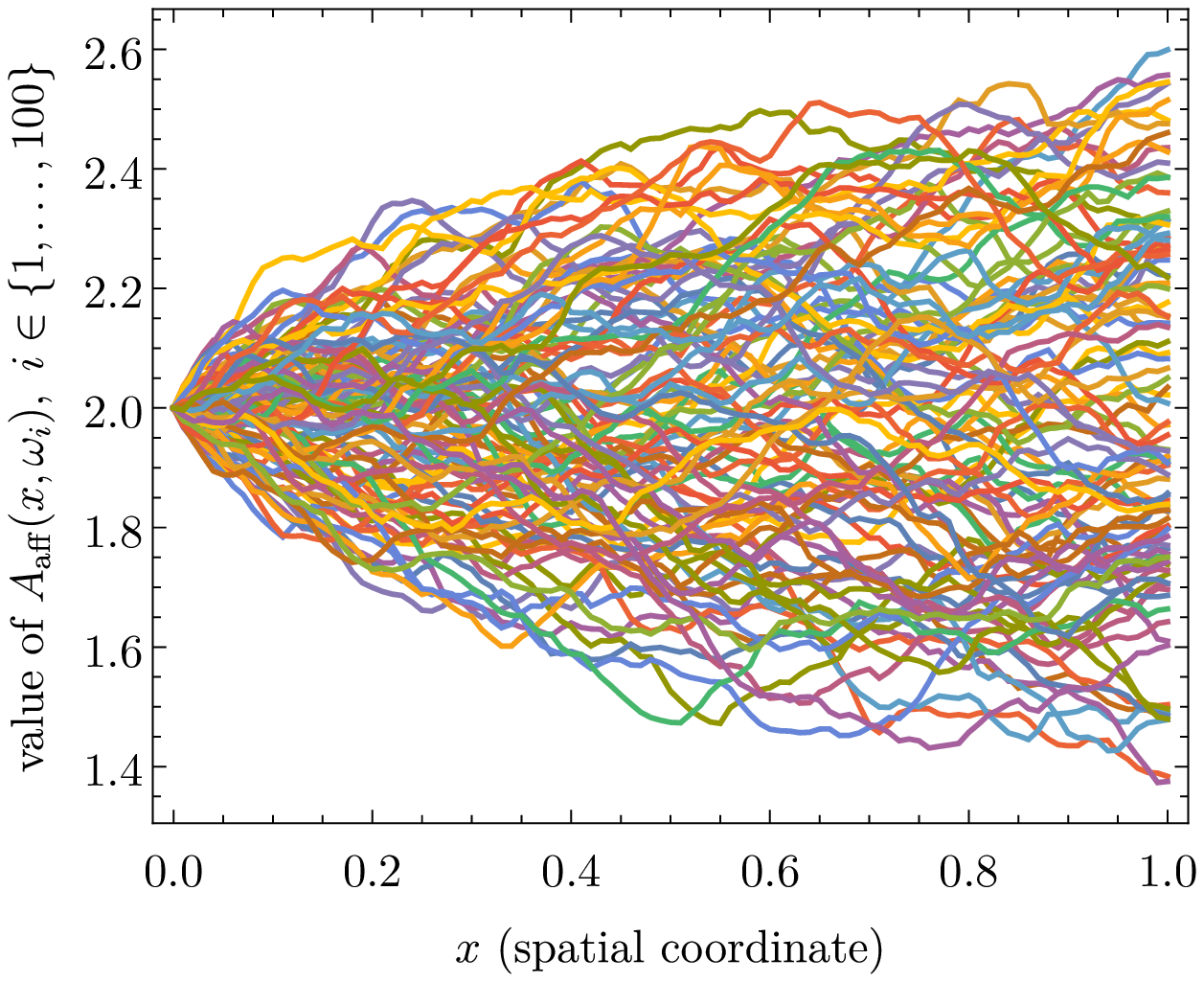}}}\quad
%{{\includegraphics[height=.36 \textwidth]{images/affine2.eps}}}\quad
\subfloat[An ensemble of realizations drawn from the periodic random field $A_{\rm per}(x,\omega)=\overline{a}(x)+\frac{1}{\sqrt{6}}\sum_{j=1}^{100}\sin(2\pi\, Y_j(\omega))\,\psi_j(x)$, $x\in[0,1{]}$, $Y_j\sim U([-\frac12,\frac12{]})$.]
{{\includegraphics[height=.36 \textwidth]{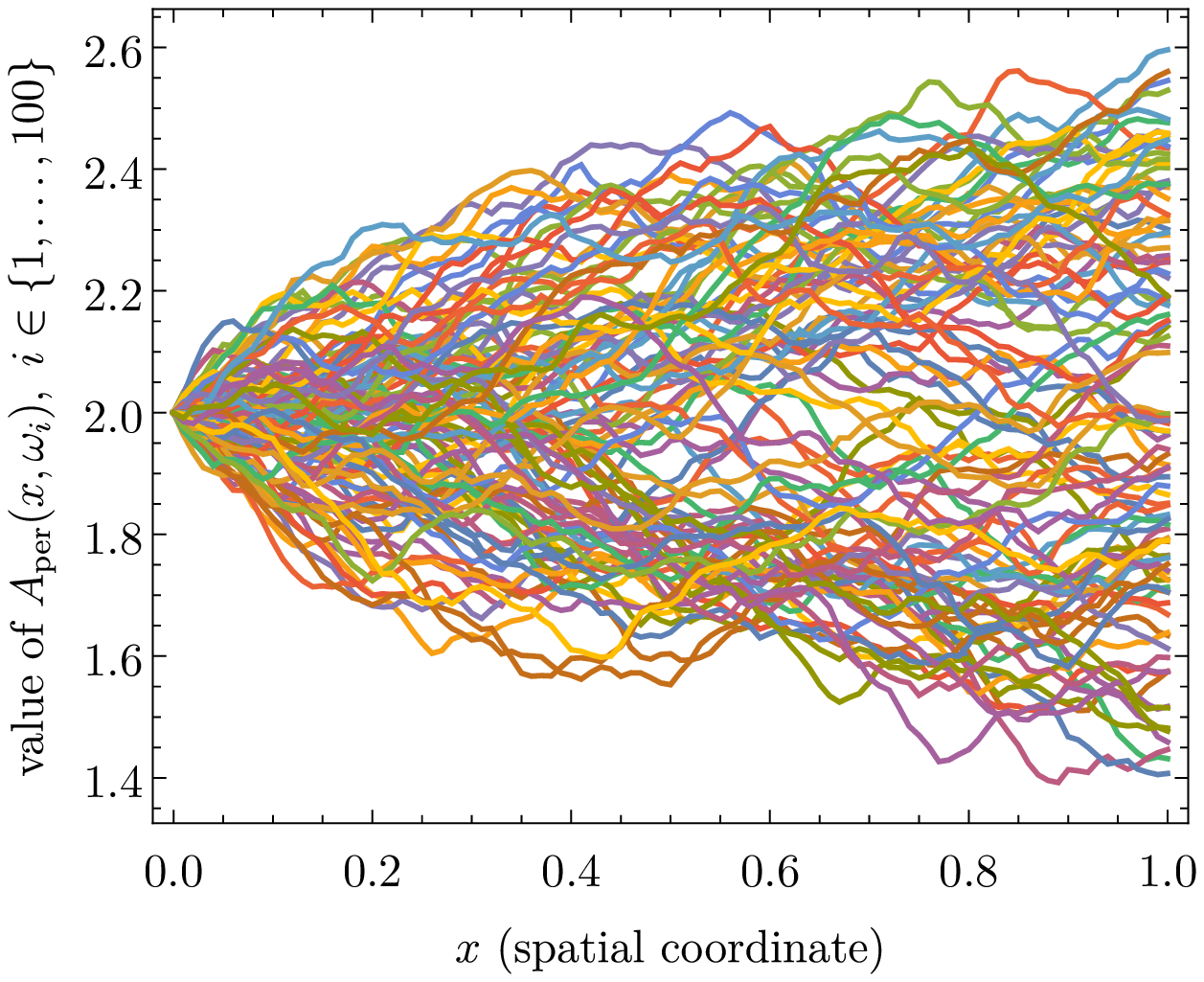}}}
%{{\includegraphics[height=.36 \textwidth]{images/periodic2.eps}}}
\caption{An illustration of $100$ individual realizations drawn from the affine and periodic random fields that correspond to the same mean $\overline{a}(x)=2$ and fluctuations $\psi_j(x)=j^{-3/2}\sin((j-\tfrac12)\pi x)$  with stochastic dimension $s=100$, which constitute a Wiener-like process in the interval $[0,1]$.}\label{fig:realizations}
\end{figure}

\begin{figure}[!t]
\centering
\subfloat
{{\includegraphics[height=.3482 \textwidth]{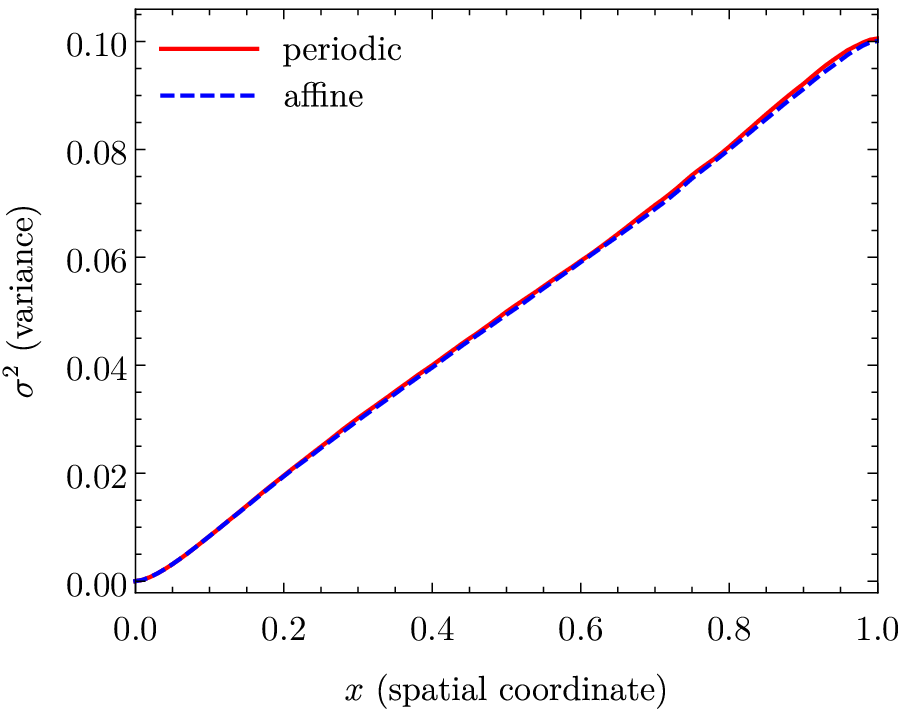}}}\quad
%{{\includegraphics[height=.347 \textwidth]{images/variance_sample.eps}}}\quad
\subfloat
{{\includegraphics[height=.3482 \textwidth]{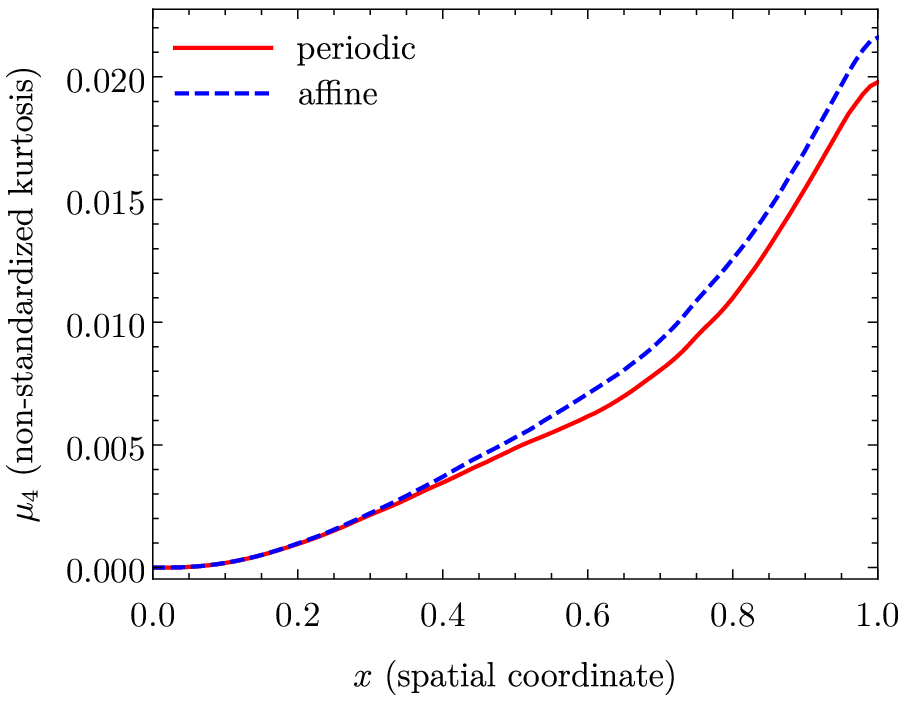}}}\quad
%{{\includegraphics[height=.347 \textwidth]{images/nonstandardized_kurtosis_sample.eps}}}\quad
\caption{A comparison of the sample variance and (nonstandardized) sample kurtosis corresponding to $10\,000$ samples drawn from the fields $A_{\rm aff}$ and $A_{\rm per}$ from Figure~\ref{fig:realizations} with $s=100$.}\label{fig:moments}
\end{figure}

An ensemble of randomly generated realizations of a pair of affine and periodic random fields is illustrated in Figure~\ref{fig:realizations}\, with the choices $\overline{a}(x):=2$ and $\psi_j(x):=j^{-3/2}\sin((j-\tfrac12)\pi x)$ for $x\in[0,1]$ and $j\in\mathbb{N}$. While the individual realizations of the affine and periodic fields are obviously different, the statistical moments of the fields coincide up to second order, but all higher moments will be different; see Figure~\ref{fig:moments}. We note that the fourth central moment (or nonstandardized kurtosis) of the random field is
\[
\mathbb{E}[(A(\boldsymbol{x},\cdot)-\overline{a}(\boldsymbol{x}))^4]=\frac{1}{96}\sum_{j\geq 1}\psi_j^4(\boldsymbol{x})+\frac{1}{24}\sum_{j\geq 1}\sum_{k\geq j+1}\psi_j^2(\bsx)\,\psi_k^2(\bsx)
\]
for the periodic model, whereas the factor $1/96$ is replaced by $1/80$ for the affine model. We are not aware of any
modeling reason to prefer the affine model over the periodic model.

To give further insight into the proposed new model, we note that we
are in effect replacing each of the countably infinite number of uniformly
distributed random variables $Y_j$ by new random variables $\sin(2\pi
Y_j)$ (when renormalized to lie between $-1$ and $+1$). It is easily seen
that the probability density function for each of these new random
variables is $1/(\pi\sqrt{1-y^2})$ for $-1\le y\le 1$, which in the
context of orthogonal polynomials is the weight function associated with
Chebyshev polynomials of the first kind.

This observation means that there is a close connection between our
new model and generalized polynomial chaos (GPC) (see
\cite{xiukarniadakis02}) in which the dependence of the solution on the
stochastic variables is expressed as a linear combination of multivariate
basis functions of orthogonal polynomials with respect to a variety of
weight functions. In the GPC setting we have merely changed from the
uniform probability distribution to the one associated with Chebyshev
polynomials of the first kind: a popular choice (see
\cite{boydbook,chkifa15,hansenschwab13,rauhutschwab17,shenyu10}) because
of the attractive simplicity of the Chebyshev polynomials. Nevertheless,
it should be emphasized that the proposed approximation scheme is
\emph{not} the same as GPC, since the input field (and hence also the
solution) is here made to be periodic, and the natural approximation is by
trigonometric rather than algebraic polynomials.

We use the random field~\eqref{eq:A_affine_sin} to model the uncertain diffusion coefficient of the following PDE problem: find $u\!:D\times \Omega\to\mathbb{R}$ that satisfies
\begin{align}\label{eq:pdeproblem}
\begin{array}{rll}
-\nabla\cdot (A(\bsx,\omega)\nabla u(\bsx,\omega))&\!\!\!\!=f(\bsx),& \bsx\in D,\\
u(\bsx,\omega)&\!\!\!\!=0,& \bsx\in \partial D,
\end{array}
\end{align}
for almost all events $\omega\in\Omega$. Then we approximate $\mathbb{E}[G(u)]$, where $G$ is a bounded, linear functional of the solution to~\eqref{eq:pdeproblem}. The motivation for the choice~\eqref{eq:A_affine_sin} of the random field is that the random field $u$ is now a
1-periodic function of the random variable $\omega$, and periodic integrands are known
to be especially advantageous in the context of so-called \emph{lattice} cubature
rules~\cite{sloanjoe94}.  By using the periodic model of random fields instead of the affine
model, it is possible to carry out lattice rule calculations of
expected values in high dimensions with higher order convergence rates, instead of being restricted, as, for example, in~\cite{kss12,kssmultilevel,schwab13}, to a convergence rate of at best $\mathcal{O}(n^{-1})$.

In this paper we leave open the choice of the fluctuations $(\psi_j)_{j\geq 1}$ but note that if the covariance function $K(\bsx,\bsx'):={\rm cov}(A)(\bsx,\bsx')$ of the field $A(\bsx,\omega)$ is specified, then the appropriate choice is to take the $\psi_j$ to be suitably normalized eigenfunctions of the integral operator with kernel $K$:
$$
\int_DK(\bsx,\bsx')\,\psi_j(\bsx')\,{\rm d}\bsx'=\lambda_j\psi_j(\bsx),\quad \bsx\in D,
$$
where $\lambda_1\geq \lambda_2\geq \cdots\geq 0$ are the eigenvalues of the integral operator, and the eigenfunctions are orthogonal with respect to the $L_2$-inner product $\langle\cdot,\cdot\rangle_{L_2}$ and normalized by
$$
\int_D|\psi_j(\bsx)|^2\,{\rm d}\bsx=12\lambda_j\quad\text{for }j\geq 1.
$$
In this case~\eqref{eq:mercerlike} becomes Mercer's theorem for the covariance function, and~\eqref{eq:A_affine_alpha} is a version of the Karhunen--Lo\`{e}ve expansion; see \cite{schwabgittelson,loeve,schwabtodor}. The only nonstandard point in the proof of the Karhunen--Lo\`{e}ve theorem is the occurrence of the function $\Theta$, but the only properties that are needed are those in~\eqref{eq:alpha2}.

The main result in this paper is as follows. We show that if in addition to~\eqref{eq:summable}, the fluctuation operators satisfy
\begin{align*}
\sum_{j\geq 1}\|\psi_j\|_{L_\infty}^p<\infty\quad\text{for some }0<p< 1,%\label{eq:psumdemo}
\end{align*}
as well as certain regularity assumptions to be made precise later, the overall error of the discretized PDE problem~\eqref{eq:pdeproblem} with the diffusion coefficient truncated to the first $s$ terms is given by
\begin{align}
\mathcal{O}(s^{-2/p+1}+h^2+n^{-1/p})\label{eq:overall}
\end{align}
using a first order finite element solver with mesh size $h$ and $n$ rank-1 lattice cubature points in $[-\tfrac12,\tfrac12]^s$ generated using the component-by-component (CBC) algorithm. The error term~\eqref{eq:overall} consists of the \emph{dimension truncation error}, \emph{first order finite element discretization error}, and the \emph{cubature error}, respectively, and the implied coefficient is {\em independent} of the truncated dimension $s$ as well as $h$ and $n$. In particular, we note that we are able to obtain a higher order lattice cubature convergence rate $\mathcal{O}(n^{-1/p})$ beating the $\mathcal{O}(n^{-\min\{1/p-1/2,1-\delta\}})$, $\delta>0$, rate for randomly shifted lattice rules in, e.g.,~\cite{kss12,schwab13}. The same rate $\mathcal{O}(n^{-1/p})$ has been obtained for the affine model in, e.g.,~\cite{spodpaper14} but with interlaced polynomial lattice rules, which are more complicated than rank-1 lattice rules in their construction. Moreover, the dimension truncation error rate $\mathcal{O}(s^{-2/p+1})$ in~\eqref{eq:overall} matches the recent result for affine-parametric operator equations~\cite{gantner}. We also discuss the case $p=1$. Higher order convergence for the finite element error can potentially be obtained by using higher order elements.

This paper is structured as follows. We present the notations and discuss the preliminaries in Subsection~\ref{sec:notations}. The periodic parametric mathematical model is introduced in Section~\ref{sec:parametricmodel}. We assess the regularity of this model with respect to the parametric variable in Subsection~\ref{sec:parametricregularity} and consider the dimension truncation error and finite element discretization errors in Subsection~\ref{sec:dimtruncfem}. The quasi-Monte Carlo (QMC) method as it applies to the periodic framework is discussed in Section~\ref{sec:periodicqmc}, and we show in Subsection~\ref{sec:higheroder} that the use of rank-1 lattice rules in the periodic setting yields a higher order convergence rate for our model provided that appropriate \emph{smoothness-driven product and order dependent} (SPOD) weights are used in the lattice rule construction. To this end, Subsection~\ref{sec:cbc} contains a description of the fast CBC algorithm for rank-1 lattice cubature rules using SPOD weights.  The overall error estimate for the discretized PDE problem is presented in Section~\ref{sec:combinederror}. We present numerical experiments in Section~\ref{sec:numex} that assess the cubature convergence rate. We end this paper with some conclusions on our results.

\subsection{Notations and preliminaries}\label{sec:notations}
We follow the convention $\mathbb{N}=\{1,2,3,\ldots\}$ and use $\mathbb{N}_0$ to denote the set of natural numbers including zero. Moreover, we use the shorthand notation $\{m:n\}=\{m,m+1,\ldots,n\}$ for integers such that $m\leq n$ and set
$$
U:=[-\tfrac12,\tfrac12]^{\mathbb{N}}\quad\text{and}\quad U_s:=[-\tfrac12,\tfrac12]^s\quad\text{for }s\in \mathbb{N}.
$$
We define the integral over the set $U$ by
$$
\int_UF(\bsy)\,{\rm d}\bsy:=\lim_{t\to\infty}\int_{U_t}F(y_1,\ldots,y_t,0,0,\ldots)\,{\rm d}y_1\cdots{\rm d}y_t.
$$
For fixed $s\in\mathbb{N}$, we introduce the set $\overline{U}_s:=\{(y_j)_{j\geq s+1}:y_j\in [-\tfrac12,\tfrac12],~j\geq s+1\}$ and write
\begin{align}
\int_{\overline{U}_s}F\big(\bsy_{\{s+1:\infty\}}\big)\,{\rm d}\bsy_{\{s+1:\infty\}}=\lim_{t\to\infty}\int_{U_t}F(y_{s+1},\ldots,y_{s+t},0,0,\ldots)\,{\rm d}y_{s+1}\cdots{\rm d}y_{s+t}\label{eq:special}
\end{align}
to mean integration over the variables $(y_j)_{j\geq s+1}$.

Let the set of all multi-indices with finite support be denoted by
\[
\mathcal{I}:=\{\boldsymbol{m}\in\mathbb{N}_0^\infty:|{\rm supp}(\boldsymbol{m})|<\infty\},
\]
where we define the support of a multi-index by ${\rm supp}(\boldsymbol{m}):=\{j\in\mathbb{N}:m_j\neq 0\}$, and $|{\rm supp}(\boldsymbol{m})|$ is the cardinality of the support. Here and throughout this manuscript, we refer to the $j{\rm th}$ component of a multi-index $\boldsymbol{m}$ as $m_j$. Moreover, we define
\[
|\boldsymbol{m}|:=\sum_{j\geq 1}m_j
\]
for multi-indices $\boldsymbol{m}\in\mathcal{I}$. Let $\boldsymbol{x}=(x_j)_{j\geq 1}$ be a sequence and $\boldsymbol{\alpha}\in\mathcal{I}$. We denote
\[
\boldsymbol{x}^{\boldsymbol{\alpha}}:=\prod_{j\in {\rm supp}(\boldsymbol{\alpha})}x_j^{\alpha_j}.
\]
In addition, we use the notation $\boldsymbol{\alpha}\leq\boldsymbol{\beta}$ to signify that $\alpha_j\leq \beta_j$ for all $j\geq 1$.

We assume in the sequel that $D\subseteq\mathbb{R}^d$, $d\in\{1,2,3\}$, is a bounded domain with a Lipschitz regular boundary. This assumption also justifies us taking the Sobolev norm of the space $H_0^1(D)$ to be
\begin{align}
\|w\|_{H_0^1}:=\|\nabla w\|_{L_2},\quad w\in H_0^1(D).\label{eq:H01def}
\end{align}
The duality pairing between $H^{-1}(D)$ of $H_0^1(D)$ is denoted by $\langle\cdot,\cdot\rangle$.

We establish the following notations and assumptions regarding the finite element approximation of $w\in H_0^1(D)$. Let us assume that $D$ is a convex and bounded polyhedron with plane faces. We denote by $\{V_h\}_h$ a family of finite element subspaces $V_h\subset H_0^1(D)$, parametrized by the mesh size $h>0$, which are spanned by continuous, piecewise linear finite element basis functions such that each $V_h$ is obtained from an initial, regular triangulation of $D$ by recursive, uniform bisection of simplices. We use the notation $w_h\in V_h$ to denote the finite element approximation of $w$ in the finite element space $V_h$.

\section{Parametric weak formulation}\label{sec:parametricmodel}
The parametric weak formulation of~\eqref{eq:pdeproblem} is, for $\bsy\in U$, to find  $u(\cdot,{\bsy})\in H_0^1(D)$ such that
\begin{align}
\int_Da(\bsx,{\bsy})\nabla u(\bsx,{\bsy})\cdot\nabla \phi(\bsx)\,{\rm d}\bsx=\langle f,\phi\rangle\quad \forall\phi\in H_0^1(D),\label{eq:weak}
\end{align}
where $f\in H^{-1}(D)$, and the diffusion coefficient is assumed to have the form
\begin{align}
a(\bsx,{\bsy})=\overline{a}(\bsx)+\frac{1}{\sqrt{6}}\sum_{j\geq 1}\sin(2\pi y_j)\psi_j(\bsx),\quad\bsx\in D,~\bsy\in U,\label{eq:rfield}
\end{align}
consistently with~\eqref{eq:A_affine_sin}. Furthermore, let $G\!:H_0^1(D)\to\mathbb{R}$ be a bounded, linear mapping. As the {\em quantity of interest}, we consider the expectation of $\boldsymbol{y}\mapsto G(u(\cdot,\boldsymbol{y}))$ taken over the parametric space:
\begin{align}
\mathbb{E}[G(u)]=\int_UG(u(\bsy))\,{\rm d}\bsy.\label{eq:expectation}%\lim_{s\to\infty}\int_{U_s}G(u(\cdot,(y_1,\ldots,y_s,0,0,\ldots)))\,{\rm d}y_1\cdots{\rm d}y_s.\label{eq:expectation}
\end{align}

We state the following assumptions which are the same as the assumptions in~\cite{kss12}:
\begin{itemize}
\item[\textnormal{(A1)}] $\overline{a}\in L_\infty(D)$ and
    $\sum_{j\geq 1}\|\psi_j\|_{L_\infty}<\infty$;
\item[\textnormal{(A2)}] there exist positive constants $a_{\rm max}$
    and $a_{\rm min}$ such that $0<a_{\rm min}\leq a(\bsx,{\bsy})\leq
    a_{\rm max}<\infty$ for all $\bsx\in D$ and $\bsy\in U$;
\item[\textnormal{(A3)}]$\sum_{j\geq 1}\|\psi_j\|_{L_\infty}^p<
    \infty$ for some $0<p\leq 1$;
\item[\textnormal{(A4)}] $\overline{a}\in W^{1,\infty}(D)$ and
    $\sum_{j\geq 1}\|\psi_j\|_{W^{1,\infty}}<\infty$, where
    $$\|v\|_{W^{1,\infty}}:=\max\{\|v\|_{L_\infty},\|\nabla
    v\|_{L_\infty}\};$$
\item[\textnormal{(A5)}] $\|\psi_1\|_{L_\infty}\geq
    \|\psi_2\|_{L_\infty}\geq\cdots$;
\item[\textnormal{(A6)}] the physical domain
    $D\subseteq\mathbb{R}^d$, $d\in\{1,2,3\}$, is a convex and bounded
    polyhedron with plane faces.
\end{itemize}
We refer to these assumptions as they are needed.

For convenience, we introduce the following notation to mean the dimensionally truncated exact solution to~\eqref{eq:weak}:
$$
u^s(\cdot,\bsy):=u(\cdot,(y_1,\ldots,y_s,0,0,\ldots))\quad\forall\bsy\in U,
$$
and we define $u_h^s(\cdot,\bsy):=u_h(\cdot,(y_1,\ldots,y_s,0,0,\ldots))\in V_h$ for all $\bsy\in U$ to mean the dimensionally truncated finite element solution to~\eqref{eq:weak}.

\subsection{Parametric regularity of the solution}\label{sec:parametricregularity}

We proceed to derive a regularity estimate for the problem~\eqref{eq:weak}
with respect to the parametric variable $\bsy$. The approach we take here
follows the argument of~\cite{kss12}, where a uniform affine model of the
uncertain diffusion coefficient was considered.

We begin by remarking that a straightforward application of the Lax--Milgram lemma ensures that~\eqref{eq:weak} is uniquely solvable over the whole parametric domain and that the solution can be bounded {a priori}.
\begin{lemma}\label{lemma:laxmilgram}
Under the assumptions \textnormal{(A1)} and \textnormal{(A2)}, the
weak formulation~\eqref{eq:weak} has a unique solution
$u(\cdot,\boldsymbol{y})\in H_0^1(D)$ for any $\boldsymbol{y}\in U$ such
that
\[
\|u(\cdot,\boldsymbol{y})\|_{H_0^1}\leq\frac{\|f\|_{H^{-1}}}{a_{\rm min}}
\]
for any source term $f\in H^{-1}(D)$.
\end{lemma}

Let $\boldsymbol{m}\in\mathcal{I}$ be a multi-index. It is easy to see that the mixed partial derivatives of~\eqref{eq:rfield} with respect to $\bsy$ are
\begin{align}
\partial^{\boldsymbol{m}}a(\bsx,\bsy)=\begin{cases}
a(\bsx,\bsy)&\text{if }\boldsymbol{m}=\mathbf{0},\\
\displaystyle\frac{1}{\sqrt{6}}(2\pi)^k\sin\!\left(2\pi y_j+k\frac{\pi}{2}\right)\psi_j(\bsx)&\text{if }\boldsymbol{m}=k\mathbf{e}_j,~k\geq 1,\\
0&\text{otherwise,}
\end{cases}\label{eq:diffcoefderiv}
\end{align}
where $\mathbf{e}_j\in\mathcal{I}$ denotes the multi-index whose $j{\rm th}$ component is $1$ and all other components are $0$. This is due to the dependence of $a$ on each $y_j$ being in separate additive terms: if we differentiate once or more with respect to $y_j$, then we obtain an expression depending only on $y_j$ and $\psi_j$, and if we then differentiate with respect to a different component of the $\bsy$ variable, we get $0$.

Let $\boldsymbol{\nu}\in\mathcal{I}$ be a multi-index with $\boldsymbol{\nu}\neq\mathbf{0}$. We differentiate~\eqref{eq:weak} on both sides to get
\begin{align*}
\int_D\partial^{\boldsymbol{\nu}}\bigg(a(\bsx,\bsy)\nabla u(\bsx,\bsy)\cdot \nabla\phi(\bsx)\bigg)\,{\rm d}\bsx=0\quad\forall\phi\in H_0^1(D),
\end{align*}
which, after an application of the Leibniz product rule, yields
\[
\int_D\bigg(\sum_{\boldsymbol{m}\leq \boldsymbol{\nu}}\binom{\boldsymbol{\nu}}{\boldsymbol{m}}(\partial^{\boldsymbol{m}}a)(\bsx,\bsy)\nabla(\partial^{\boldsymbol{\nu}-\boldsymbol{m}}u(\bsx,\bsy))\cdot \nabla\phi(\bsx)\bigg){\rm d}\bsx=0\quad\forall\phi\in H_0^1(D).
\]
Plugging in~\eqref{eq:diffcoefderiv} and separating out the case $\boldsymbol{m}=\mathbf{0}$, we obtain
\begin{align*}
&\int_D a(\bsx,\bsy)\nabla(\partial^{\boldsymbol{\nu}}u(\bsx,\bsy))\cdot \nabla \phi(\bsx)\,{\rm d}\bsx\\
&=-\sum_{j\geq 1}\sum_{k=1}^{\nu_j}\int_D\binom{\nu_j}{k}\frac{(2\pi)^k}{\sqrt{6}}\sin\!\big(2\pi y_j+k\frac{\pi}{2}\big)\psi_j(\bsx)\nabla( \partial^{\boldsymbol{\nu}-k\mathbf{e}_j}u(\bsx,\bsy))\cdot\nabla \phi(\bsx)\,{\rm d}\bsx
\end{align*}
for all $\phi\in H_0^1(D)$. In particular, we can choose to test this formula  against $\phi=(\partial^{\boldsymbol{\nu}}u)(\cdot,\bsy)$. By applying the ellipticity assumption $a(\bsx,\bsy)\geq a_{\rm min}$ on the left-hand side and $|\psi_j(\bsx)|\leq \|\psi_j\|_{L_\infty}$ as well as the Cauchy--Schwarz inequality on the right-hand side, we obtain
\begin{align*}
&a_{\rm min}\int_D|\nabla (\partial^{\boldsymbol{\nu}}u)(\bsx,\bsy)|^2\,{\rm d}\bsx\\
&\leq\! \sum_{j\geq 1}\!\sum_{k=1}^{\nu_j}\!\binom{\nu_j}{k}\!\frac{(2\pi)^k}{\sqrt{6}}\|\psi_j\|_{L_\infty}\!\bigg(\int_D|\nabla\partial^{\boldsymbol{\nu}-k\mathbf{e}_j}u(\bsx,\bsy)|^2\,{\rm d}\bsx\bigg)^{\!\!1/2}\!\!\bigg(\int_D|\nabla\partial^{\boldsymbol{\nu}}u(\bsx,\bsy)|^2\,{\rm d}\bsx\bigg)^{\!\!1/2}\!.
\end{align*}
Eliminating the common factor on both sides and using~\eqref{eq:H01def} yields for $\bsnu\neq\mathbf{0}$
\begin{align}
\|\partial^{\boldsymbol{\nu}}u(\cdot,\bsy)\|_{H_0^1}\leq \sum_{j\geq 1}\sum_{k=1}^{\nu_j}\binom{\nu_j}{k}(2\pi)^kb_j\| \partial^{\boldsymbol{\nu}-k\mathbf{e}_j}u(\cdot,\bsy)\|_{H_0^1},\label{eq:recurrence}
\end{align}
where we set
\begin{align}
b_j:=\frac{1}{\sqrt{6}}\frac{\|\psi_j\|_{L_\infty}}{a_{\rm min}}\quad\text{for }j\in\mathbb{N}.\label{eq:bseq}
\end{align}
This differs from the definition of $b_j$ in \cite{kss12} by the
factor $1/\sqrt{6}$.

 Our goal is to use the recurrence~\eqref{eq:recurrence} to derive an
explicit upper bound on the term
$\|\partial^{\boldsymbol{\nu}}u(\cdot,\bsy)\|_{H_0^1}$ for all
$\boldsymbol{\nu}\in\mathcal{I}$. It turns out that \emph{Stirling numbers
of the second kind} (or \emph{Stirling partition numbers}) play a large
role in the forthcoming analysis; they are defined by
\[
S(n,k):=\frac{1}{k!}\sum_{j=0}^k(-1)^{k-j}\binom{k}{j}j^n
\]
for $n\geq k\geq 0$ except for $S(0,0):=1$. The bound on $\|\partial^{\boldsymbol{\nu}}u(\cdot,\bsy)\|_{H_0^1}$ (see Theorem~\ref{lemma:inductivebound} below) follows from the following result, stated in a general form in case it is useful in other contexts.
\begin{lemma}\label{lemma:inductivebound0}
Let $B,c>0$, and let $(\mathbb{A}_{\boldsymbol{\nu}})_{\boldsymbol{\nu}\in\mathcal{I}}$ and $(\Upsilon_j)_{j\in\mathbb{N}}$ be sequences of non-negative real numbers that satisfy the recurrence
\begin{align}
\mathbb{A}_{\mathbf{0}}\leq B\quad\text{and}\quad \mathbb{A}_{\boldsymbol{\nu}}\leq \sum_{j\geq 1}\sum_{k=1}^{\nu_j}\binom{\nu_j}{k}c^k\Upsilon_j\mathbb{A}_{\boldsymbol{\nu}-k\textbf{e}_j}\quad\text{for }\boldsymbol{\nu}\in{\mathcal{I}}\setminus\{\mathbf{0}\}.\label{eq:recu2}
\end{align}
Then
\begin{align}
\mathbb{A}_{\boldsymbol{\nu}}\leq c^{|\boldsymbol{\nu}|}B\sum_{\boldsymbol{m}\leq\boldsymbol{\nu}}|\boldsymbol{m}|!\,\boldsymbol{\Upsilon}^{\boldsymbol{m}}\prod_{i\geq 1}S(\nu_i,m_i).\label{eq:induc2}
\end{align}
Moreover, if equalities hold in the formulae~\eqref{eq:recu2}, then there is equality in~\eqref{eq:induc2}.
\end{lemma}
\proof We prove this result by carrying out an induction argument on $|\boldsymbol{\nu}|$ based on the recurrence~\eqref{eq:recu2}. The base step $\boldsymbol{\nu}=\mathbf{0}$ is resolved immediately. For arbitrary $\boldsymbol{\nu}\in\mathcal{I}\setminus\{\boldsymbol{0}\}$, suppose that the claim holds for all multi-indices of order $<|\boldsymbol{\nu}|$. In particular, if $\nu_j\geq k\geq 1$ for some $j\geq 1$, then the induction hypothesis gives
\[
\mathbb{A}_{\boldsymbol{\nu}-k\textbf{e}_j}\leq c^{|\boldsymbol{\nu}|-k}\,B\sum_{\boldsymbol{m}\leq\boldsymbol{\nu}-k\textbf{e}_j}|\boldsymbol{m}|!\,\boldsymbol{\Upsilon}^{\boldsymbol{m}}S(\nu_j-k,m_j)\prod_{\substack{i\geq 1\\ i\neq j}}S(\nu_i,m_i).
\]
Applying the recursion~\eqref{eq:recu2} in conjunction with the inequality above yields
\begin{align} \label{eq:step1}
\mathbb{A}_{\boldsymbol{\nu}}\leq
c^{|\boldsymbol{\nu}|}B\sum_{j\geq 1}\sum_{k=1}^{\nu_j}\binom{\nu_j}{k}\Upsilon_j\sum_{\boldsymbol{m}\leq\boldsymbol{\nu}-k\textbf{e}_j}|\boldsymbol{m}|!\,\boldsymbol{\Upsilon}^{\boldsymbol{m}}S(\nu_j-k,m_j)\prod_{\substack{i\geq 1\\ i\neq j}}S(\nu_i,m_i).
\end{align}

For given $\boldsymbol{m}$, $\bsnu$, $\boldsymbol{\Upsilon}$ and an
index $j$, we define
$\boldsymbol{m'}=(m_1,\ldots,m_{j-1},m_{j+1},\ldots)$,
$\boldsymbol{\nu'}=(\nu_1,\ldots,\nu_{j-1},\nu_{j+1},\ldots)$, and
$\boldsymbol{\Upsilon'}=(\Upsilon_1,\ldots,\Upsilon_{j-1},\Upsilon_{j+1},\ldots)$,
respectively. Then we may write the $j{\rm th}$ term in the outer sum from
\eqref{eq:step1} as
\begin{align} \label{eq:step2}
& \sum_{k=1}^{\nu_j}\binom{\nu_j}{k}\sum_{m_j=0}^{\nu_j-k}\Upsilon_j^{m_j+1}S(\nu_j-k,m_j)\sum_{\boldsymbol{m'}\leq\boldsymbol{\nu'}}(|\boldsymbol{m'}|+m_j)!\,\boldsymbol{\Upsilon'}^{\boldsymbol{m'}}\prod_{\substack{i\geq 1\\ i\neq j}}S(\nu_i,m_i) \nonumber\\
&= \sum_{\boldsymbol{m'}\leq\boldsymbol{\nu'}}\boldsymbol{\Upsilon'}^{\boldsymbol{m'}}
   \bigg(\prod_{\substack{i\geq 1\\ i\neq j}}S(\nu_i,m_i)\bigg) \sum_{k=1}^{\nu_j}\binom{\nu_j}{k}\sum_{m_j=0}^{\nu_j-k}\Upsilon_j^{m_j+1}S(\nu_j-k,m_j)(|\boldsymbol{m'}|+m_j)! \nonumber\\
&= \sum_{\boldsymbol{m'}\leq\boldsymbol{\nu'}}\boldsymbol{\Upsilon'}^{\boldsymbol{m'}}
   \bigg(\prod_{\substack{i\geq 1\\ i\neq j}}S(\nu_i,m_i)\bigg)\sum_{m_j=0}^{\nu_j-1}\Upsilon_j^{m_j+1}(|\boldsymbol{m'}|+m_j)!\sum_{k=1}^{\nu_j-m_j}\binom{\nu_j}{k}S(\nu_j-k,m_j),
\end{align}%
where we swapped the order of the sums over $k$ and $m_j$. Furthermore, it holds that
$$
\sum_{k=1}^{\nu_j-m_j}\binom{\nu_j}{k}S(\nu_j-k,m_j)=(m_j+1)S(\nu_j,m_j+1)\quad\text{for}~m_j<\nu_j,
$$
which can be verified either by direct calculation based on the definition
of $S(n,k)$ or as a consequence of~\cite[equation~(9.25)]{stirling}.
Thus \eqref{eq:step2} becomes
\begin{align*}
\sum_{\boldsymbol{m'}\leq\boldsymbol{\nu'}}\boldsymbol{\Upsilon'}^{\boldsymbol{m'}}\prod_{\substack{i\geq 1\\ i\neq j}}S(\nu_i,m_i)\sum_{m_j=0}^{\nu_j-1}\Upsilon_j^{m_j+1}(|\boldsymbol{m'}|+m_j)!\,(m_j+1)S(\nu_j,m_j+1),
\end{align*}
and together with \eqref{eq:step1} this yields
\begin{align*}
\mathbb{A}_{\boldsymbol{\nu}}
\leq c^{|\boldsymbol{\nu}|}B\sum_{j\geq 1}\sum_{\boldsymbol{m}\leq \boldsymbol{\nu}-\mathbf{e}_j}|\boldsymbol{m}|!\,(m_j+1)\Upsilon_j\boldsymbol{\Upsilon}^{\boldsymbol{m}} S(\nu_j,m_j+1)\prod_{\substack{i\geq 1\\ i\neq j}}S(\nu_i,m_i).
\end{align*}
Since $S(k,0)=0$ for all $k\geq 1$, a straightforward computation shows that
\begin{align*}
\sum_{j\geq 1}\sum_{\boldsymbol{m}\leq \boldsymbol{\nu}-\mathbf{e}_j}|\boldsymbol{m}|!\,(m_j+1)\Upsilon_j\boldsymbol{\Upsilon}^{\boldsymbol{m}}S(\nu_j,m_j+1)\prod_{\substack{i\geq 1\\ i\neq j}}S(\nu_i,m_i)\\
=\sum_{\boldsymbol{m}\leq\boldsymbol{\nu}}|\boldsymbol{m}|!\,\boldsymbol{\Upsilon}^{\boldsymbol{m}}\prod_{i\geq 1}S(\nu_i,m_i),
\end{align*}
which simplifies the upper bound into
\[
\mathbb{A}_{\boldsymbol{\nu}}\leq c^{|\boldsymbol{\nu}|}B\sum_{\boldsymbol{m}\leq\boldsymbol{\nu}}|\boldsymbol{m}|!\,\boldsymbol{\Upsilon}^{\boldsymbol{m}}\prod_{i\geq 1}S(\nu_i,m_i),
\]
completing the proof.\quad\endproof

The desired result can be obtained as an immediate corollary to Lemma~\ref{lemma:inductivebound0} using Lemma~\ref{lemma:inductivebound0} and~\eqref{eq:recurrence}.

\begin{theorem}\label{lemma:inductivebound}
Under the assumptions \textnormal{(A1)} and \textnormal{(A2)}, for any
$\bsy\in U$, let $u(\cdot,\bsy)\in H_0^1(D)$ be the solution of the
problem~\eqref{eq:weak} with the source term $f\in H^{-1}(D)$, and let
$\boldsymbol{b}=(b_j)_{j\geq 1}$ be the sequence defined
by~\eqref{eq:bseq}. Then for any multi-index
$\boldsymbol{\nu}\in\mathcal{I}$ we have
\[
\|\partial^{\boldsymbol{\nu}}u(\cdot,\bsy)\|_{H_0^1}\leq\frac{\|f\|_{H^{-1}}}{a_{\rm min}}(2\pi)^{|\boldsymbol{\nu}|} \sum_{\boldsymbol{m}\leq\boldsymbol{\nu}}|\boldsymbol{m}|!\,\boldsymbol{b}^{\boldsymbol{m}}\prod_{i\geq 1}S(\nu_i,m_i).
\]
The result also holds for the dimension-truncated finite element solution $u_h^s(\cdot,\bsy)\in V_h$ for all $s\in \mathbb{N}$, $\bsy\in U$.
\end{theorem}
\subsection{Dimension truncation and finite element discretization errors}\label{sec:dimtruncfem}
In practice, it is generally only possible to solve the problem~\eqref{eq:weak} approximately using, e.g., the finite element method and with the series~\eqref{eq:rfield} truncated to finitely many terms. In this section, we discuss the approximation errors caused by the finite element discretization and dimension truncation.

In the affine setting, the fundamental dimension truncation error bound
has already been discussed  in~\cite{kss12} leading to an error bound of
the order $\mathcal{O}(s^{-2/p+2})$. While this analysis can also be
applied to the periodic setting with only minuscule changes to the
argument, Gantner~\cite{gantner} has recently proved an improved bound of
order $\mathcal{O}(s^{-2/p+1})$ in the context of affine-parametric
operator equations. In the following, we prove an analogous result for the
problem~\eqref{eq:weak}--\eqref{eq:expectation}. While the proof technique
we use is the same as in~\cite{gantner}, we present the proof for
completeness in order to highlight that the result holds also in the
periodic framework. The following proof also differs from~\cite{gantner}
insofar as we do \emph{not} need to put a restriction on the size of the
sum, e.g., $\sum_{j\geq 1}b_j<\sqrt{6}$.

\begin{lemma}[cf.~{\cite[Theorem~1]{gantner}}]\label{lemma:dimtrunc}
Under the assumptions \textnormal{(A1)--(A3)} and
\textnormal{(A5)}, for any $\bsy\in U$, let $u(\cdot,\bsy)\in H_0^1(D)$
denote the solution to the problem~\eqref{eq:weak} with the source term
$f\in H^{-1}(D)$, and let $G\in H^{-1}(D)$. If $0<p<1$, then for any
$s\in\mathbb{N}$ there exists a constant $C>0$ such that
\[
\bigg|\int_UG(u(\cdot,\boldsymbol{y})-u^s(\cdot,\boldsymbol{y}))\,{\rm d}\boldsymbol{y}\bigg|\leq C\|G\|_{H^{-1}}\|f\|_{H^{-1}}s^{-2/p+1}.
\]
If $p=1$, then
$$
\bigg|\int_UG(u(\cdot,\bsy)-u^s(\cdot,\bsy))\,{\rm d}\bsy\bigg|\leq C\|G\|_{H^{-1}}\|f\|_{H^{-1}}\bigg(\sum_{j\geq s+1}b_j\bigg)^2.
$$
In both cases, $C>0$ denotes a generic constant that does not depend on $s$, $f$, or $G$.
\end{lemma}

\proof We define the operators $A(\bsy)\!:H_0^1(D)\to H^{-1}(D)$
for $\bsy\in U$, and $A_j\!:H_0^1(D)\to H^{-1}(D)$ for
$j\in\mathbb{N}$, by setting $\langle A(\bsy)w,\phi\rangle := \langle
a(\cdot,\bsy)\nabla w,\nabla\phi\rangle_{L_2}$ for all $\phi\in
H_0^1(D)$, and $\langle A_j w,\phi\rangle := \langle
\frac{1}{\sqrt{6}}\psi_j\nabla w,\nabla \phi\rangle_{L_2}$ for all
$\phi\in H_0^1(D)$, respectively. Moreover, we define
$A^s(\bsy):=A((y_1,\ldots,y_s,0,0,\ldots))$ and denote
$u(\bsy):=u(\cdot,\bsy)$ and $u^s(\bsy):=u^s(\cdot,\bsy)$ for all $s\in
\mathbb{N}$, $\bsy\in U$. These definitions lead to the identity
\begin{align*}
A(\bsy)-A^s(\bsy)=\sum_{j\geq s+1}\sin(2\pi y_j)A_j\quad\forall\bsy\in U,~s\in\mathbb{N}.
\end{align*}
Let $w\in H_0^1(D)$. Lemma~\ref{lemma:laxmilgram} and~\eqref{eq:H01def}
together with
$$
\|A(\bsy)w\|_{H^{-1}}=\sup_{\phi\in H_0^1(D)\setminus\{0\}}\frac{\langle a(\cdot,\bsy)\nabla w,\nabla \phi\rangle_{L_2}}{\|\phi\|_{H_0^1}}\leq a_{\rm max}\|w\|_{H_0^1}
$$
imply that both operators $A(\bsy)$ and $A^s(\bsy)$ are boundedly
invertible linear maps for all $\bsy\in U$. Furthermore, we obtain
\begin{align*}
\|A^s(\bsy)^{-1}A_jw\|_{H_0^1}\leq \frac{\|A_jw\|_{H^{-1}}}{a_{\rm min}}=\frac{1}{a_{\rm min}}\sup_{\phi\in H_0^1(D)\setminus\{0\}}\frac{\langle\frac{1}{\sqrt{6}}\psi_j\nabla w,\nabla \phi\rangle_{L_2}}{\|\phi\|_{H_0^1}}\leq b_j\|w\|_{H_0^1},
\end{align*}
where the sequence $(b_j)_{j\geq 1}$ is defined as in~\eqref{eq:bseq}.
In consequence, this yields
\begin{align}
&\sup_{\bsy\in U}\|A^s(\bsy)^{-1}A_j\|_{\mathscr{L}(H_0^1(D))}\leq b_j,\label{eq:futureproof}\\
&\sup_{\bsy\in U}\|A^s(\bsy)^{-1}(A(\bsy)-A^s(\bsy))\|_{\mathscr{L}(H_0^1(D))}\leq \sum_{j\geq s+1}b_j.\label{eq:neumannconv}
\end{align}
In what follows, we omit the argument $\bsy$ and denote the operator norm by
$\|\cdot\|=\|\cdot\|_{\mathscr{L}(H_0^1(D))}$ for brevity.

Since the sequence $(b_j)_{j\geq 1}$ is summable, there exists
$s'\in\mathbb{N}$ such that for all $s\geq s'$ the upper bound in
\eqref{eq:neumannconv} is at most $1/2$. Let us assume that $s\ge
s'$. For future reference, we note that this implies for all $s\geq
s'$
\begin{align}
b_j \le \frac12\quad\forall j\geq s+1\quad\text{and}\quad
\sum_{j\geq s+1}b_j^2
\le \sum_{j\geq s+1}b_j
\le \frac12.\label{eq:twoinequs}
\end{align}%

It follows from~\eqref{eq:neumannconv} and our assumption $s\geq s'$ that the Neumann series
\begin{align*}
A^{-1}=(I+(A^s)^{-1}(A-A^s))^{-1}(A^s)^{-1}=\sum_{k\geq 0}(-(A^s)^{-1}(A-A^s))^k(A^s)^{-1}
\end{align*}
is well defined. Moreover, we have the representation
\begin{align}
\int_UG(u\!-\!u^s)\,{\rm d}\bsy
&=\int_U G((A^{-1}\!-\!(A^s)^{-1})f)\,{\rm d}\bsy=\sum_{k\geq 1}\int_UG((-(A^s)^{-1}(A\!-\!A^s))^ku^s)\,{\rm d}\bsy\notag\\
&=\sum_{k\geq 1}(-1)^k\int_UG\bigg(\bigg(\sum_{j\geq s+1}\sin(2\pi y_j)(A^s)^{-1}A_j\bigg)^ku^s\bigg)\,{\rm d}\bsy.\label{eq:split}
\end{align}

For each $k\in\mathbb{N}$, we note that the integrand in~\eqref{eq:split}
can be expanded as
\begin{align*}
\bigg(\sum_{j\geq s+1}\sin(2\pi y_j)(A^s)^{-1}A_j\bigg)^k=\sum_{\eta_1,\ldots,\eta_k\geq s+1}\bigg(\prod_{i=1}^k\sin(2\pi y_{\eta_i})\bigg)\bigg(\prod_{i=1}^k(A^s)^{-1}A_{\eta_i}\bigg),
\end{align*}
where the product symbol is assumed to respect the order of the
noncommutative operators. Using the independence of the components of $\bsy\in
U$ and~\eqref{eq:special}, the integral over $U$ in~\eqref{eq:split} can be written as a product of integrals
\begin{align*}
&\int_U G\bigg(\bigg(\sum_{j\geq s+1}\sin(2\pi y_j)(A^s)^{-1}A_j\bigg)^ku^s\bigg)\,{\rm d}\bsy\\
%&= \int_U \sum_{\eta_1,\ldots,\eta_k\geq s+1}\!\!\! G \bigg(
%   \bigg(\prod_{i=1}^k\sin(2\pi y_{\eta_i})\bigg)
%   \bigg(\prod_{i=1}^k(A^s)^{-1}A_{\eta_i}\bigg)u^s\bigg)\,{\rm d}\bsy\\
&= \!\!\sum_{\eta_1,\ldots,\eta_k\geq s+1}\!\!\!
   \bigg(\underbrace{\int_{\overline{U}_s} \prod_{i=1}^k\sin(2\pi y_{\eta_i})\,{\rm d}\bsy_{\{s+1:\infty\}}}_{=: I_1}\bigg)
   \!\bigg(\underbrace{\int_{U_s}\!G\bigg(\bigg(\prod_{i=1}^k(A^s)^{-1}A_{\eta_i}\bigg)u^s\bigg)\,{\rm d}\bsy_{\{1:s\}}}_{=: I_2}\bigg),
\end{align*}
where $I_1\ge 0$ because it can be written as a product of univariate
integrals of the form $\int_{-1/2}^{1/2} \sin(2\pi y_j)^m \,\rd y_j$ for
$m\in\bbN$, which take values between $0$ and $1$ (importantly, this expression is
zero when $m=1$), while we can estimate $I_2$ by
\begin{align*}
  |I_2|
\le \|G\|_{H^{-1}} \bigg(\prod_{i=1}^k\sup_{\bsy\in U_s}\|(A^s)^{-1}A_{\eta_i}\|\bigg)\|u^s\|_{H_0^1}
\leq \frac{\|G\|_{H^{-1}}\|f\|_{H^{-1}}}{a_{\rm min}} \bigg(\prod_{i=1}^k b_{\eta_i}\bigg).
\end{align*}
Thus
\begin{align*}
& \bigg| (-1)^k \int_UG\bigg(\bigg(\sum_{j\geq s+1}\sin(2\pi y_j)(A^s)^{-1}A_j\bigg)^ku^s\bigg)\,{\rm d}\bsy\bigg|\\
&\leq \frac{\|G\|_{H^{-1}}\|f\|_{H^{-1}}}{a_{\rm min}}\sum_{\eta_1,\ldots,\eta_k\geq s+1}
 \bigg(\int_{\overline{U}_s} \prod_{i=1}^k\sin(2\pi y_{\eta_i})\,{\rm d}\bsy_{\{s+1:\infty\}}\bigg)
 \bigg(\prod_{i=1}^k b_{\eta_i}\bigg)\\
&= \frac{\|G\|_{H^{-1}}\|f\|_{H^{-1}}}{a_{\rm min}} \int_{\overline{U}_s}
 \sum_{\eta_1,\ldots,\eta_k\geq s+1}
 \bigg(\prod_{i=1}^k\sin(2\pi y_{\eta_i})\bigg)
 \bigg(\prod_{i=1}^k b_{\eta_i}\bigg) \,{\rm d}\bsy_{\{s+1:\infty\}}\\
&=\frac{\|G\|_{H^{-1}}\|f\|_{H^{-1}}}{a_{\rm min}}\int_{\overline{U}_s}\bigg(\sum_{j\geq s+1}\sin(2\pi y_j)b_j\bigg)^k\,{\rm d}\bsy_{\{s+1:\infty\}}\\
&=\frac{\|G\|_{H^{-1}}\|f\|_{H^{-1}}}{a_{\rm min}} \int_{\overline{U}_s}
  \sum_{\substack{|\boldsymbol{\nu}|=k\\ \nu_j=0~\forall j\leq s}}
   \frac{k!}{\boldsymbol{\nu}!}\,
   \bigg(\prod_{j\geq s+1} \sin(2\pi y_j)^{\nu_j}\bigg)\,
   \bigg(\prod_{j\geq s+1} b_j^{\nu_j}\bigg) \,{\rm d}\bsy_{\{s+1:\infty\}} \\
%&= \frac{\|G\|_{H^{-1}}\|f\|_{H^{-1}}}{a_{\rm min}}\!\!\sum_{\substack{|\boldsymbol{\nu}|=k\\ \nu_j=0~\forall j\leq s}}
%   \frac{k!}{\boldsymbol{\nu}!}\,
%   \bigg(\prod_{j\geq s+1} \int_{-1/2}^{1/2} \sin(2\pi y_j)^{\nu_j}\,{\rm d} y_j\bigg)\,\boldsymbol{b}^{\bsnu}\\
&\leq\frac{\|G\|_{H^{-1}}\|f\|_{H^{-1}}}{a_{\rm min}}\!\!\sum_{\substack{|\boldsymbol{\nu}|=k\\ \nu_j=0~\forall j\leq s\\ \nu_j\neq 1~\forall j\geq 1}}
  \frac{k!}{\boldsymbol{\nu}!}\,\boldsymbol{b}^{\bsnu},
\end{align*}
where we have used the multinomial theorem together with
$\boldsymbol{\nu}!:=\prod_{i\geq 1}\nu_i!$ for
$\boldsymbol{\nu}\in\mathcal{I}$, Lemma~\ref{lemma:laxmilgram}, and the
bound~\eqref{eq:futureproof}. The key observation is that this term vanishes whenever any component of $\boldsymbol{\nu}$ is equal to $1$, and consequently the
term vanishes when $k=1$.

We may now estimate~\eqref{eq:split} by splitting the sum into the
$k\ge k'$ terms and the $k < k'$ terms for a value of $k'$ to be specified
later. We obtain
\begin{align*}
\bigg|\int_UG(u-u^s)\,{\rm d}\bsy\bigg|
\leq \frac{\|G\|_{H^{-1}}\|f\|_{H^{-1}}}{a_{\rm min}}\bigg(
 \sum_{k\geq k'}\bigg(\sum_{j\geq s+1}b_j\bigg)^k
 + k'!\sum_{2\le k<k'} \sum_{\substack{|\boldsymbol{\nu}|=k\\ \nu_j=0~\forall j\leq s\\ \nu_j\neq 1~\forall j\geq 1}}\boldsymbol{b}^{\bsnu}
 \bigg).
\end{align*}

Consider first the case $0<p<1$. The $k\ge k'$ terms can be
bounded using the geometric series as
$$
\sum_{k\geq k'}\bigg(\sum_{j\geq s+1}b_j\bigg)^k\leq \bigg(\sum_{j\geq s+1}b_j\bigg)^{k'}\frac{1}{1-\sum_{j\geq s+1}b_j}
\leq C_1s^{k'(-1/p+1)},
$$
where we used the inequality $\sum_{j\geq s+1}b_j\leq (\sum_{j\geq
1}b_j^p)^{1/p}s^{-1/p+1}$ (see~\cite[Theorem~5.1]{kss12}), and the ensuing constant $C_1:=2(\sum_{j\geq
1}b_j^p)^{k'/p}$ is independent of $s$, $f$, and $G$. On the other hand,
for each $2\le k<k'$ we use the estimate
\begin{align*}
\sum_{\substack{|\boldsymbol{\nu}|=k\\ \nu_j=0~\forall j\leq s\\ \nu_j\neq 1~\forall j\geq 1}}\boldsymbol{b}^{\boldsymbol{\nu}}
&\leq \sum_{\substack{0\ne|\boldsymbol{\nu}|_\infty\leq k\\ \nu_j=0~\forall j\leq s\\ \nu_j\neq 1~\forall j\geq 1}}\boldsymbol{b}^{\boldsymbol{\nu}}=\prod_{j\geq s+1}\bigg(1+\sum_{\ell=2}^kb_j^\ell\bigg)-1=\prod_{j\geq s+1}\bigg(1+b_j^2\frac{1-b_j^{k-1}}{1-b_j}\bigg)-1\notag\\
&\leq \prod_{j\geq s+1}\big(1+2 b_j^2\big)-1\leq \exp\bigg(2\sum_{j\geq s+1}b_j^2\bigg)-1\leq C_2s^{-2/p+1},
\end{align*}
where we used both inequalities in~\eqref{eq:twoinequs}, the inequalities
${\rm e}^x\leq 1+({\rm e}-1)x$ for all $x\in[0,1]$ and $\sum_{j\geq
s+1}b_j^2\leq\frac{1}{2/p-1}(\sum_{j\geq 1}b_j^p)^{1/p}s^{-2/p+1}$, and
the resulting constant $C_2:=\frac{2({\rm e}-1)}{2/p-1}(\sum_{j\geq
1}b_j^p)^{1/p}$ is independent of $s$, $f$, and $G$. Hence we conclude
that
\begin{align*}
\bigg|\int_UG(u-u^s)\,{\rm d}\bsy\bigg|
\leq \frac{\|G\|_{H^{-1}}\|f\|_{H^{-1}}}{a_{\rm min}}\Big(
  C_1s^{k'(-1/p+1)} \,+\, k'!\,(k'-2)C_2s^{-2/p+1} \Big).
\end{align*}
We therefore choose $k' := \lceil(2-p)/(1-p)\rceil$ to balance the two
terms. This proves the assertion for $s\geq s'$ after a trivial
adjustment of the constant factors. The result can be extended to all
$s\in\mathbb{N}$ by noticing that
\begin{align*}
\bigg|\int_UG(u-u^s)\,{\rm d}\boldsymbol{y}\bigg|\leq 2\frac{\|G\|_{H^{-1}}\|f\|_{H^{-1}}}{a_{\rm min}}\leq 2\frac{\|G\|_{H^{-1}}\|f\|_{H^{-1}}}{a_{\rm min}(s'-1)^{-2/p+1}}s^{-2/p+1}
\end{align*}
holds for all $1\leq s<s'$, and the claim follows by a trivial adjustment of all of the constants involved.

For $p=1$ we amend the above argument slightly to obtain
$$
\bigg|\int_UG(u(\bsy)-u^s(\bsy))\,{\rm d}\bsy\bigg|\leq C\|G\|_{H^{-1}}\|f\|_{H^{-1}}\bigg(\sum_{j\geq s+1}b_j\bigg)^2,
$$
where $C>0$ is a constant independent of $s$, $f$, and $G$.\quad\endproof

Regarding the finite element approximation error, it is clear that an analogous result to the one presented in~\cite{kss12} holds.

\begin{lemma}[cf.~{\cite[Theorem~5.1]{kss12}}]\label{lemma:femerror}
Under assumptions \textnormal{(A1)}, \textnormal{(A2)}, \textnormal{(A4)},
and \textnormal{(A6)}, for any $\bsy\in U$, let $u(\cdot,\bsy)\in
H_0^1(D)$ denote the solution to~\eqref{eq:weak} with the source term
$f\in H^{-1+t}(D)$ such that $0\leq t\leq 1$, and let $G\in H^{-1+t'}(D)$
with $0\leq t'\leq 1$. Then the finite element approximations satisfy the following
asymptotic convergence estimate as $h\to 0$:
\[
|G(u(\cdot,\boldsymbol{y})-u_h(\cdot,\boldsymbol{y}))|\leq Ch^{t+t'}\|f\|_{H^{-1+t}}\|G\|_{H^{-1+t'}},
\]
where $0\leq t+t'\leq 2$ and the constant $C>0$ is independent of $h$ and $\boldsymbol{y}$.
\end{lemma}

\emph{Remark}. We note that the limiting case $t=t'=1$ in Lemma~\ref{lemma:femerror} corresponds to taking $f\in L_2(D)$ and $G\in L_2(D)$, where the dual of $L_2(D)$ is identified with itself, resulting in a convergence rate of $\mathcal{O}(h^2)$.

\section{QMC in the periodic setting}\label{sec:periodicqmc}
QMC methods are a class of numerical methods designed  to approximate multivariate integrals such as
\[
I_s(F)=\int_{[0,1]^s}F(\bsy)\,{\rm d}\bsy
\]
for a continuous integrand $F$ by using an equal weight cubature formula of the form
\[
Q_{s,n}(F)=\frac{1}{n}\sum_{k=0}^{n-1}F(\boldsymbol{y}_{k}),
\]
where $\boldsymbol{y}_{0},\ldots,\boldsymbol{y}_{n-1}\in[0,1]^s$ are prescribed cubature nodes.

We consider rank-1 lattice rules, where the QMC nodes $\Lambda:=\{\boldsymbol{y}_{0},\ldots,\boldsymbol{y}_{n-1}\}$ are taken to be of the form
\[
\boldsymbol{y}_{k}=\left\{\frac{k\boldsymbol{z}}{n}\right\},\quad k\in\{0,\ldots,n-1\},
\]
where $\{\boldsymbol{x}\}$ denotes taking the componentwise fractional part of $\boldsymbol{x}\in\mathbb{R}^s$ and $\boldsymbol{z}\in\mathbb{N}^s$ is called the {\em generating vector} of a lattice rule. It is well known that the lattice rule error for functions with absolutely convergent Fourier series is precisely~\cite{sk87}
\begin{align}
Q_{s,n}(F)-I_s(F)=\sum_{\boldsymbol{h}\in\Lambda^\perp\setminus\{\mathbf{0}\}}\hat F(\boldsymbol{h}),\label{eq:laterr}
\end{align}
where $\hat F(\boldsymbol{h}):=\int_{[0,1]^s}F(\bsy){\rm e}^{-2\pi{\rm i}\bsy\cdot\boldsymbol{h}}\,{\rm d}\bsy$ for $\boldsymbol{h}\in\mathbb{Z}^s$ and we denote the \emph{dual lattice} by $\Lambda^\perp=\Lambda^\perp(\boldsymbol{z})=\{\boldsymbol{h}\in\mathbb{Z}^s: \boldsymbol{h}\cdot\boldsymbol{z}\equiv 0~({\rm mod}~n)\}$, which is defined with respect to the generating vector $\boldsymbol{z}$ of the rank-1 lattice rule.

Let $F\in C([0,1)^s)$ be a 1-periodic function with respect to each of its variables, and set
$$
r_\alpha(\bsgamma,\boldsymbol{h}):=\gamma_{{\rm supp}(\boldsymbol{h})}^{-1}\prod_{j\in{\rm supp}(\boldsymbol{h})}|h_j|^\alpha\quad\text{for}~\alpha>1~\text{and}~\boldsymbol{h}\in\mathbb{Z}^s,
$$
where ${\rm supp}(\boldsymbol{h}):=\{j\in\{1:s\}:h_j\neq 0\}$ and
$\bsgamma=(\gamma_{\mathfrak{u}})_{\mathfrak{u}\subseteq\{1:s\}}$ denotes
a collection of nonnegative weights. Using the error
formula~\eqref{eq:laterr}, we obtain
\begin{align}
|I_s(F)-Q_{s,n}(F)|=\bigg|\sum_{\boldsymbol{h}\in\Lambda^\perp\setminus\{\mathbf{0}\}}\hat F(\boldsymbol{h})\frac{r_\alpha(\bsgamma,\boldsymbol{h})}{r_\alpha(\bsgamma,\boldsymbol{h})}\bigg|\leq P_\alpha(\bsgamma,\boldsymbol{z})\|F\|_{\alpha},\label{eq:holderbound}
\end{align}
where the factor depending only on the QMC nodes is defined by
\[
P_{\alpha}(\bsgamma,\boldsymbol{z}):=\sum_{\boldsymbol{h}\in \Lambda^\perp\setminus\{\mathbf{0}\}}\frac{1}{r_\alpha(\bsgamma,\boldsymbol{h})}\quad\text{for }\alpha>1,
\]
and the norm is given by
$$\|F\|_\alpha:=\sup_{\boldsymbol{h}\in\mathbb{Z}^s}|\hat F(\boldsymbol{h})|r_\alpha(\bsgamma,\boldsymbol{h})\quad\text{for }\alpha>1.$$
Since the inequality~\eqref{eq:holderbound} is sharp, we see that
$P_\alpha(\bsgamma,\boldsymbol{z})$  is the worst-case error in the
space with $\|F\|_\alpha\leq 1$. The quantity
$P_\alpha(\bsgamma,\boldsymbol{z})$ is well known in classical
lattice rule literature (at least for the unweighted case
$\gamma_{\mathfrak{u}}\equiv 1$, see~\cite{sloanjoe94}) and coincides
with the squared error term in the Hilbert space setting considered in the
paper~\cite{korobovpaper06}, leading us to conclude the following.

\begin{lemma}\label{lemma:cbcerrorbound}
Let $s\in\mathbb{N}$ and prime $n$, and let $\bsgamma=(\gamma_{\mathfrak{u}})_{\mathfrak{u}\subseteq\{1:s\}}$ be a collection of nonnegative weights. Let $F\in C([0,1)^s)$ be a {\rm 1}-periodic function with respect to each of its variables such that $\|F\|_{\alpha}<\infty$. Then a generating vector $\boldsymbol{z}\in\mathbb{N}^s$ can be constructed by the CBC algorithm such that
\begin{align*}
|I_s(F)-Q_{s,n}(F)|\leq\bigg(\frac{1}{n-1}\sum_{\varnothing\neq\mathfrak{u}\subseteq\{1:s\}}\gamma_{\mathfrak{u}}^\lambda (2\zeta(\alpha\lambda))^{|\mathfrak{u}|}\bigg)^{1/\lambda}\|F\|_\alpha
\end{align*}
for $\lambda\in(1/\alpha,1]$. Here, $\zeta(x):=\sum_{k\geq 1}k^{-x}$
denotes the Riemann zeta function for $x>1$.
\end{lemma}

\proof It can be readily verified that $F$ has an absolutely convergent
Fourier series given that $\|F\|_\alpha<\infty$. The claim then follows
from the previous discussion in conjunction
with~\cite[Theorem~5]{korobovpaper06}.\quad\endproof

The result can be extended to nonprime $n$ by replacing $n-1$ with
Euler's totient function $\varphi_{\rm
tot}(n):=|\{m\in\mathbb\{1:n-1\}:{\rm gcd}(m,n)=1\}|$. In particular,
$1/\varphi_{\rm tot}(n) \leq 2/n$ if $n$ is a prime power.

When $\alpha\geq 2$ is an integer, it can be shown that
\begin{align}
\|F\|_{\alpha}\leq \max_{\mathfrak{u}\subseteq\{1:s\}}\frac{1}{(2\pi)^{\alpha|\mathfrak{u}|}}\frac{1}{\gamma_{\mathfrak{u}}}\int_{[0,1]^{|\mathfrak{u}|}}\bigg|\int_{[0,1]^{s-|\mathfrak{u}|}}\bigg(\prod_{j\in\mathfrak{u}}\frac{\partial}{\partial y_j}\bigg)^{\alpha}F(\boldsymbol{y})\,{\rm d}\boldsymbol{y}_{\{1:s\}\setminus\mathfrak{u}}\bigg|{\rm d}\boldsymbol{y}_{\mathfrak{u}}\label{eq:normupperbound}
\end{align}
provided that $F$ has mixed partial derivatives of order $\alpha$. Furthermore, when $\alpha$ is even, we can write
\begin{align}
P_\alpha(\bsgamma,\boldsymbol{z})=\frac{1}{n}\sum_{k=0}^{n-1}\sum_{\varnothing\neq \mathfrak{u}\subseteq\{1:s\}}\gamma_{\mathfrak{u}}\prod_{j\in\mathfrak{u}}\omega\bigg(\bigg\{\frac{kz_j}{n}\bigg\}\bigg),\label{eq:palfa}
\end{align}
where
$$
\omega(x):=(2\pi)^\alpha\frac{B_\alpha(x)}{(-1)^{\alpha/2+1}\alpha!}\quad\text{for }x\in[0,1],
$$
and $B_\alpha$ denotes the Bernoulli polynomial of degree $\alpha$.

\subsection{Higher order convergence in the PDE context}\label{sec:higheroder} In this section, we let the assumptions \textnormal{(A1)}--\textnormal{(A3)} be in effect. We are interested in the expectation of the functional $F(\boldsymbol{y}):=G(u(\cdot,\boldsymbol{y}-\tfrac{\mathbf{1}}{\mathbf{2}}))$, where $G$ denotes a bounded, linear functional $G\!:H_0^1(D)\to\mathbb{R}$, $u(\cdot,\boldsymbol{y}-\tfrac{\mathbf{1}}{\mathbf{2}})\in H_0^1(D)$ is the solution to the weak formulation~\eqref{eq:weak}, and we let $\boldsymbol{y}\in [0,1]^{\mathbb{N}}$.

For an integer $\alpha\geq 2$, we estimate the norm as follows:
\begin{align*}
\bigg|\bigg(\prod_{j\in \mathfrak{u}}\frac{\partial}{\partial y_j}\bigg)^\alpha F(\bsy)\bigg|
&=\bigg|\bigg(\prod_{j\in\mathfrak{u}}\frac{\partial}{\partial y_j}\bigg)^{\alpha}G(u(\cdot,\boldsymbol{y}-\tfrac{\mathbf{1}}{\mathbf{2}}))\bigg|=\bigg|G\bigg(\bigg(\prod_{j\in\mathfrak{u}}\frac{\partial}{\partial y_j}\bigg)^{\alpha}u(\cdot,\boldsymbol{y}-\tfrac{\mathbf{1}}{\mathbf{2}})\bigg)\bigg|\\
&\leq \|G\|_{H^{-1}}\left\|\bigg(\prod_{j\in\mathfrak{u}}\frac{\partial}{\partial y_j}\bigg)^{\alpha}u(\cdot,\boldsymbol{y}-\tfrac{\mathbf{1}}{\mathbf{2}})\right\|_{H_0^1}.
\end{align*}
We thus obtain, using~\eqref{eq:normupperbound} and Theorem~\ref{lemma:inductivebound},
\begin{align*}
\|F\|_{\alpha}\leq \frac{\|G\|_{H^{-1}}\|f\|_{H^{-1}}}{a_{\rm min}}\max_{\mathfrak{u}\subseteq\{1:s\}}\frac{1}{\gamma_{\mathfrak{u}}}\sum_{\boldsymbol{m}_{\mathfrak{u}}\in\{1:\alpha\}^{|\mathfrak{u}|}}|\boldsymbol{m}_{\mathfrak{u}}|!\prod_{j\in \mathfrak{u}}(b_j^{m_j}S(\alpha,m_j))
\end{align*}
since $S(\alpha,0)=0$ for $\alpha\neq 0$.  We now choose the weights to be
\begin{align} \label{eq:weights}
\gamma_{\mathfrak{u}}=\sum_{\boldsymbol{m}_{\mathfrak{u}}\in\{1:\alpha\}^{|\mathfrak{u}|}}|\boldsymbol{m}_{\mathfrak{u}}|!\prod_{j\in\mathfrak{u}}\big(b_j^{m_j}S(\alpha,m_j)\big)\quad\forall\mathfrak{u}\subseteq\{1:s\},
\end{align}
which ensures that $\|F\|_{\alpha}$ is bounded. These weights have a very
specific form: they are  \emph{SPOD weights}, first seen
in~\cite{spodpaper14}. We then observe that  the bound for the error term
in Lemma~\ref{lemma:cbcerrorbound} becomes
\begin{align*}
|I_s(F)-Q_{s,n}(F)|\leq\frac{\|G\|_{H^{-1}}\|f\|_{H^{-1}}}{a_{\rm min}} \bigg(\frac{2}{n}\bigg)^{1/\lambda}C(s,\alpha,\lambda),
\end{align*}
where
\[
C(s,\alpha,\lambda):=\bigg(\sum_{\varnothing\neq\mathfrak{u}\subseteq\{1:s\}}\bigg(\sum_{\boldsymbol{m}_{\mathfrak{u}}\in\{1:\alpha\}^{|\mathfrak{u}|}}|\boldsymbol{m}_{\mathfrak{u}}|!\prod_{j\in\mathfrak{u}}\big(b_j^{m_j}S(\alpha,m_j)\big)\bigg)^\lambda (2\zeta(\alpha\lambda))^{|\mathfrak{u}|}\Bigg)^{1/\lambda}
\]
for $\lambda\in(1/\alpha,1]$ and $n$  a prime power.

Finally, we need to choose $\lambda$ in such a way that $C(s,\alpha,\lambda)$ is bounded \emph{independently} of~$s$. First applying the inequality (cf.~\cite[Theorem~19]{jenseninequ}) $$\sum_ka_k\leq \bigg(\sum_ka_k^\lambda\bigg)^{1/\lambda},\quad 0<\lambda\leq 1,~a_k\geq 0,$$  to the inner sum of $C(s,\alpha,\lambda)$ and denoting $S_{\max}(\alpha):=\max_{k\in\{1:\alpha\}}S(\alpha,k)$ yields
\begin{align}
&[C(s,\alpha,\lambda)]^\lambda\leq \sum_{\varnothing\neq\mathfrak{u}\subseteq\{1:s\}}\sum_{\boldsymbol{m}_{\mathfrak{u}}\in\{1:\alpha\}^{|\mathfrak{u}|}}(|\boldsymbol{m}_{\mathfrak{u}}|!)^\lambda\prod_{j\in\mathfrak{u}}\big(b_j^{m_j}S(\alpha,m_j)\big)^\lambda (2\zeta(\alpha\lambda))^{|\mathfrak{u}|}\notag\\
&\leq \sum_{\varnothing\neq\mathfrak{u}\subseteq\{1:s\}}\sum_{\boldsymbol{m}_{\mathfrak{u}}\in\{1:\alpha\}^{|\mathfrak{u}|}}(|\boldsymbol{m}_{\mathfrak{u}}|!)^\lambda\prod_{j\in\mathfrak{u}}\big(S_{\max}(\alpha)(2\zeta(\alpha\lambda))^{1/\lambda}b_j^{m_j}\big)^{\lambda}\notag\\
&\leq \sum_{\varnothing\neq\mathfrak{u}\subseteq\{1:s\}}\sum_{\boldsymbol{m}_{\mathfrak{u}}\in\{1:\alpha\}^{|\mathfrak{u}|}}  \bigg(|\boldsymbol{m}_{\mathfrak{u}}|!\prod_{j\in\mathfrak{u}}\beta_j^{m_j}\bigg)^{\lambda},\notag
\end{align}
where we have set $\beta_j:=\max\{1,S_{\max}(\alpha)(2\zeta(\alpha\lambda))^{1/\lambda}\}b_j$. We recast the double sum as a sum over multi-indices $\bsnu$:
\[
[C(s,\alpha,\lambda)]^\lambda\leq\sum_{\mathbf{0}\neq\boldsymbol{\nu}\in\{0:\alpha\}^s}\bigg(|\boldsymbol{\nu}|!\prod_{\substack{j=1}}^s\beta_j^{\nu_j}\bigg)^\lambda.
\]

Let us define the sequence $d_j=\beta_{\lceil j/\alpha\rceil}$, $j\geq 1$. In concrete terms, this means that
\begin{align*}
%&d_1=d_2=\cdots=d_{\alpha}=\beta_1,\\
%&d_{\alpha+1}=d_{\alpha+2}=\cdots=d_{2\alpha}=\beta_2,\\
%&\vdots\\
d_{k\alpha+1}=\,d_{k\alpha+2}=\cdots=d_{(k+1)\alpha}=\beta_{{k}+1},\quad k\in\mathbb{N}_0.
\end{align*}
We relate this definition to $C(s,\alpha,\lambda)$ by observing that
\begin{align*}
\sum_{\mathbf{0}\neq\boldsymbol{\nu}\in\{0:\alpha\}^s}\bigg(|\boldsymbol{\nu}|!\prod_{\substack{j=1}}^s\beta_j^{\nu_j}\bigg)^\lambda\leq  \sum_{\substack{\mathfrak{v}\subseteq\mathbb{Z}_+\\ |\mathfrak{v}|<\infty}}\bigg(|\mathfrak{v}|!\prod_{j\in\mathfrak{v}}d_j\bigg)^\lambda=\sum_{\ell\geq 0}(\ell!)^\lambda\sum_{\substack{\mathfrak{v}\subseteq\mathbb{Z}_+\\ |\mathfrak{v}|=\ell}}\prod_{j\in\mathfrak{v}}d_j^\lambda
\\
\leq\sum_{\ell\geq 0}(\ell!)^{\lambda-1}\bigg(\sum_{j\geq 1}d_j^\lambda\bigg)^{\ell}.
\end{align*}
The final inequality holds because $(\sum_{j\geq 1}d_j^\lambda)^\ell$
includes all the products of the form
$\prod_{j\in\mathfrak{v}}d_j^\lambda$ with $|\mathfrak{v}|=\ell$, but
since the order in which the terms in the product appear does not matter,
we can divide this by $\ell!$.

We now choose $\lambda=p$ and verify that the last expression is
finite with this choice of $\lambda$. Our assumption that
$(\|\psi_j\|_{L_\infty})_{j\geq 1}\in\ell^p$ for some $p\in(0,1]$
implies that $\sum_{j\geq 1} b_j^p<\infty$. For the inner sum, we now have
\begin{align*}
T:=\sum_{j\geq 1} d_j^p=\alpha\sum_{j\geq 1} \beta_j^p=\alpha(\max\{1,S_{\max}(\alpha)(2\zeta(\alpha p))^{1/p}\})^p\sum_{j\geq 1} b_j^p<\infty,
\end{align*}
provided that $\alpha p>1$. If additionally $p<1$,
then the ratio test implies convergence of the outer sum since
\[
\frac{((\ell+1)!)^{p-1}}{(\ell!)^{p-1}}\frac{T^{\ell+1}}{T^{\ell}}=(\ell+1)^{p-1}T\xrightarrow{\ell\to\infty}0.
\]
If $p=1$, then the sum is geometric and converges if and only if
$T<1$. An equivalent condition is that
\begin{align}
\sum_{j\geq 1}\|\psi_j\|_{L_\infty}<\frac{\sqrt{6}\,a_{\rm min}}{2\alpha\zeta(\alpha) \max_{k\in\{1:\alpha\}}S(\alpha,k)}\quad\text{for integer }\alpha\geq 2.\label{eq:pone}
\end{align}

Since the condition $1/\alpha<p\leq 1$ needs to be in effect, we conclude
that by choosing $\alpha := \lfloor 1/p\rfloor+1$ we obtain
$\mathcal{O}(n^{-1/p})$ convergence with an implied constant independent
of $s$. If $p=1$, then we assume additionally that~\eqref{eq:pone} holds.

\subsection{CBC construction with SPOD weights}\label{sec:cbc}

We describe the CBC construction of lattice rules for SPOD weights of the
general form
\begin{equation*}% \label{eq:POD}
  \gamma_\setu =
  \sum_{\boldsymbol{m}_\setu\in \{1:\sigma\}^{|\setu|}}
  \Gamma_{|\boldsymbol{m}_\setu|}
  \prod_{j\in\setu} \gamma_{j,m_j}\,,%\quad\text{for all }\mathfrak{u}\subseteq\{1:s\},
\end{equation*}
which is specified by a \emph{smoothness degree} $\sigma \in \bbN$, a
sequence $(\Gamma_\ell)_{\ell\ge 0}$, plus a sequence
$(\gamma_{j,m})_{j\ge1}$ for every $m=1,\ldots,\sigma$. Note that for
$\setu = \varnothing$, we use the convention that the empty product is
one, and we interpret the sum over $\boldsymbol{m}_\varnothing$ as a sum
with a single term $\bszero$, so that $\gamma_\varnothing = \Gamma_0$
(which in turn is typically set to $1$).

The choice of weights \eqref{eq:weights} corresponds to the specific case
$\sigma = \alpha$, $\Gamma_\ell=\ell!$, and $\gamma_{j,m_j} = b_j^{m_j}
S(\sigma,m_j)$. We consider a generic search criterion $P(\bsz)$ which
takes the same form as \eqref{eq:palfa} but with a generic function
$\omega\!:[0,1]\to \mathbb{R}$. Substituting in the weights, we can write
\begin{align*}
  P(\boldsymbol{z})
  &=
  \frac{1}{n} \sum_{k=0}^{n-1}
  \sum_{\varnothing\ne\setu\subseteq\{1:s\}}
  \sum_{\boldsymbol{m}_\setu\in \{1:\sigma\}^{|\setu|}}
  \Gamma_{|\boldsymbol{m}_\setu |} \prod_{j\in\setu} \bigg( \gamma_{j,m_j}\,
  \omega\bigg(\bigg\{\frac{kz_j}{n}\bigg\}\bigg)\bigg) \\
  &=
  \frac{1}{n} \sum_{k=0}^{n-1}
  \sum_{\satop{\bsnu\in \{0:\sigma\}^s}{\bsnu\ne\bszero}}
   \Gamma_{|\bsnu |} \prod_{\satop{j=1}{\nu_j\ne 0}}^s \bigg( \gamma_{j,\nu_j}\,
  \omega\bigg(\bigg\{\frac{kz_j}{n}\bigg\}\bigg)\bigg) \\
  &=
  \frac{1}{n} \sum_{k=0}^{n-1} \sum_{\ell=1}^{\sigma s}
  \underbrace{
  \sum_{\satop{\bsnu\in \{0:\sigma\}^s}{|\bsnu|=\ell}}
  \Gamma_{\ell} \prod_{\satop{j=1}{\nu_j\ne 0}}^s \bigg( \gamma_{j,\nu_j}\,
  \omega\bigg(\bigg\{\frac{kz_j}{n}\bigg\}\bigg)\bigg)
  }_{=:\,p_{s,\ell}(k)}
  \,.
\end{align*}
Next we find a recursive definition for $p_{s,\ell}(k)$. By considering
whether or not $\nu_s$ is zero, we can write
\begin{align} \label{eq:recur}
  p_{s,\ell}(k)
  &=
  \, \sum_{\nu_s=0}^{\min\{\ell,\sigma\}}
  \sum_{\satop{\bsnu\in \{0:\sigma\}^{s-1}}{|\bsnu|=\ell-\nu_s}}
  \Gamma_{\ell}  \prod_{\satop{j=1}{\nu_j\ne 0}}^s \bigg(\gamma_{j,\nu_j}\,
  \omega\bigg(\bigg\{\frac{kz_j}{n}\bigg\}\bigg)\bigg) \nonumber\\
  &=
  \sum_{\satop{\bsnu\in \{0:\sigma\}^{s-1}}{|\bsnu|=\ell}}
 \Gamma_{\ell} \prod_{\satop{j=1}{\nu_j\ne 0}}^{s-1} \bigg( \gamma_{j,\nu_j}\,
  \omega\bigg(\bigg\{\frac{kz_j}{n}\bigg\}\bigg)\bigg) \nonumber\\
  &\quad + \sum_{\nu_s = 1}^{\min\{\ell,\sigma\}}
  \sum_{\satop{\bsnu\in \{0:\sigma\}^{s-1}}{|\bsnu|=\ell-\nu_s}}
  \Gamma_{\ell}\,
  \gamma_{s,\nu_s}\,
  \omega\bigg(\left\{\frac{kz_s}{n}\right\}\bigg)
  \prod_{\satop{j=1}{\nu_j\ne 0}}^{s-1} \bigg(\gamma_{j,\nu_j}\,
  \omega\bigg(\bigg\{\frac{kz_j}{n}\bigg\}\bigg)\bigg) \nonumber\\
  &=
  \sum_{\satop{\bsnu\in \{0:\sigma\}^{s-1}}{|\bsnu|=\ell}}
   \Gamma_{\ell} \prod_{\satop{j=1}{\nu_j\ne 0}}^{s-1} \bigg( \gamma_{j,\nu_j}\,
  \omega\bigg(\bigg\{\frac{kz_j}{n}\bigg\}\bigg)\bigg) \nonumber\\
  &\quad + \omega\bigg(\bigg\{\frac{kz_s}{n}\bigg\}\bigg)
  \sum_{w = 1}^{\min\{\ell,\sigma\}}
  \frac{ \Gamma_{\ell}}{\Gamma_{\ell-w}}\,
  \gamma_{s,w}
  \sum_{\satop{\bsnu\in \{0:\sigma\}^{s-1}}{|\bsnu|=\ell-w}}
  \Gamma_{\ell-w}
  \prod_{\satop{j=1}{\nu_j\ne 0}}^{s-1} \bigg(\gamma_{j,\nu_j}\,
  \omega\bigg(\bigg\{\frac{kz_j}{n}\bigg\}\bigg)\bigg) \nonumber\\
  &=
  p_{s-1,\ell}(k)
  + \omega\bigg(\left\{\frac{kz_s}{n}\right\}\bigg)
  \sum_{w = 1}^{\min\{\ell,\sigma\}}
  \frac{\Gamma_{\ell}}{\Gamma_{\ell-w}}\,
  \gamma_{s,w}\, p_{s-1,\ell-w}(k).
\end{align}
Thus we have
\begin{align} \label{eq:error}
  & P(\boldsymbol{z})
  = \frac{1}{n} \sum_{k=0}^{n-1} \sum_{\ell=1}^{\sigma s}
  \bigg(
  p_{s-1,\ell}(k)
  + \omega\bigg(\bigg\{\frac{kz_s}{n}\bigg\}\bigg)
  \sum_{w = 1}^{\min\{\ell,\sigma\}}
  \frac{\Gamma_{\ell}}{\Gamma_{\ell-w}}\,
  \gamma_{s,w}\, p_{s-1,\ell-w}(k)
  \bigg) \nonumber\\
  &=
  \frac{1}{n} \sum_{k=0}^{n-1} \sum_{\ell=1}^{\sigma s}
  p_{s-1,\ell}(k)
  \!+\! \frac{1}{n} \sum_{k=0}^{n-1} \omega\bigg(\bigg\{\frac{kz_s}{n}\bigg\}\bigg)
  \sum_{\ell=1}^{\sigma s} \sum_{w = 1}^{\min\{\ell,\sigma\}}
 \! \frac{\Gamma_{\ell}}{\Gamma_{\ell-w}}
  \gamma_{s,w}\, p_{s-1,\ell-w}(k).
\end{align}
Note that the first term in~\eqref{eq:recur} is exactly the value of
$P(z_1,\ldots,z_{s-1})$ in the first $s-1$ dimensions, but this is
irrelevant for the construction.

Let $\bbZ_n:=\{0,1,\ldots,n-1\}$ denote the set of the integers modulo
$n$, and let $\mathbb{U}_n := \{ u \in \bbZ_n : \gcd(u,n) = 1 \}$ denote the
multiplicative group of integers modulo $n$ with $|\mathbb{U}_n| = \varphi_{\rm
tot}(n)$. We define the matrix
\begin{align} \label{eq:Omega}
        \bsOmega_n
&:= \left[
        \omega\bigg(\bigg\{ \frac{kz}{n} \bigg\} \bigg)
        \right]_{\substack{z \in \mathbb{U}_n \\ k \in \bbZ_n}}
  = \left[
        \omega\bigg( \frac{kz \bmod{n}}{n} \bigg)
        \right]_{\substack{z \in \mathbb{U}_n \\ k \in \bbZ_n}}
\end{align}
and the vectors
\begin{align} \label{eq:vectors}
  \bsp_{s,\ell} := \left[ p_{s,\ell}(k)\right]_{k\in\bbZ_N},
  \qquad \ell=1,\ldots,\sigma s\,,
\end{align}
where the entries $p_{s,\ell}(k)$ are defined recursively by
\eqref{eq:recur} together with $p_{s,0}(k):= 1$ for all $k$.

At step $s$, we see from \eqref{eq:error} that the CBC algorithm should
pick the value of $z_s\in \mathbb{U}_n$ which corresponds to the smallest entry in
the matrix-vector product
\[
  \bsOmega_n\,\bsx\,,
  \quad\mbox{with}\quad \bsx
  \,:=\,
  \sum_{\ell=1}^{\sigma s} \sum_{w = 1}^{\min\{\ell,\sigma\}}
  \frac{\Gamma_{\ell}}{\Gamma_{\ell-w}}\,
  \gamma_{s,w}\, \bsp_{s-1,\ell-w}\,.
\]
Then it is clear from \eqref{eq:recur} that the vectors $\bsp_{s,\ell}$
for the next iteration can be obtained recursively via
\[
        \bsp_{s,\ell}
:= \bsp_{s-1,\ell} + \bsOmega_n(z_s) \,.\!*\,
  \bigg(
  \sum_{w = 1}^{\min\{\ell,\sigma\}}
  \frac{\Gamma_{\ell}}{\Gamma_{\ell-w}}\,
  \gamma_{s,w}\, \bsp_{s-1,\ell-w}
  \bigg)\,,
\]
where $\bsOmega_n(z_s)$ denotes the row of $\bsOmega_n$ corresponding to
the chosen $z_s$, and the operator $.*$ denotes the elementwise vector
multiplication. Since the vectors $\bsp_{s-1,\ell}$ are no longer needed
in the next iteration, we can simply overwrite $\bsp_{s-1,\ell}$ with
$\bsp_{s,\ell}$. Hence, starting with the vectors $\bsp_{0,\ell}:=
\bsone_n$ requires
$
  \calO(\sigma s\, n)$ storage overall.

The fast implementation is based on ordering the indices $z \in
\mathbb{U}_n$ and $k \in \bbZ_n$ in~\eqref{eq:Omega}
and~\eqref{eq:vectors} to allow fast
matrix-vector multiplication using FFT; see~\cite{CKN06,NC06a,NC06b} for details. % The
%cost for each step of the construction requires the following
%\[
%\begin{tabular}{|l|l|}
%  \hline
%  obtain the vector $\bsx$ & $\calO(\fk{\sigma}^2 s\,n)$ operations \\
%  \hline
%  fast matrix-vector multiplication $\bsOmega_n\,\bsx$ & $\calO(n\log n)$ operations \\
%  \hline
%  update the vectors $\bsp_{s,\ell}$ & $\calO(\fk{\sigma}^2 s\, n)$ operations\\
%  \hline
%\end{tabular}
%\]
 The overall CBC construction cost is
$
  \calO(s\,n\log n + \sigma^2 s^2\,n)$ operations.

%Note that for higher order polynomial lattice rules, we expect the $N$ in
%the cost bound to become $N^\alpha$. It may be possible to reduce this to
%$\alpha N$.

\section{Combined error analysis}\label{sec:combinederror}
The overall error of the PDE problem~\eqref{eq:weak} is a combination of the \emph{dimension truncation error}, \emph{finite element discretization error}, and \emph{QMC cubature error} as
\begin{align*}
&\bigg|\int_UG(u(\cdot,\boldsymbol{y}))\,{\rm d}\bsy-\frac{1}{n}\sum_{i=0}^{n-1}G(u^s(\cdot,\boldsymbol{y}_{i}-\tfrac{\boldsymbol{1}}{\boldsymbol{2}}))\bigg|\\
&\leq \bigg|\int_UG(u(\cdot,\bsy)-u^s(\cdot,\bsy))\,{\rm d}\bsy\bigg|+\bigg|\int_{U_s}G(u^s(\cdot,\bsy)-u_h^s(\cdot,\bsy))\,{\rm d}\bsy\bigg|\\
&\quad+\bigg|\int_{U_s}G(u_h^s(\cdot,\bsy))\,{\rm d}\bsy-\frac{1}{n}\sum_{i=0}^{n-1}G(u_h^s(\cdot,\bsy_{i}-\tfrac{\boldsymbol{1}}{\boldsymbol{2}}))\bigg|,
\end{align*}
where $(\bsy_i)_{i=0}^{n-1}$ are QMC nodes in $[0,1]^s$, $u$ denotes the solution to~\eqref{eq:weak}, $u^s$ and $u_h^s$ denote the dimension-truncated solution and the corresponding finite element solution, and $G\!:H_0^1(D)\to\mathbb{R}$ is a bounded, linear functional.

We can combine the results of the previous sections to produce the following overall error bound.
\begin{theorem}
For any $\bsy\in U$, let $u(\cdot,\bsy)\in H_0^1(D)$ denote the solution
to~\eqref{eq:weak} with the source term $f\in H^{-1+t}(D)$ for some $0\leq
t\leq 1$, and let $G\in H^{-1+t'}(D)$ for some $0\leq t'\leq 1$. Let
$(\bsy_k)_{k=0}^{n-1}$ be the lattice cubature nodes in $[0,1]^s$
generated by the CBC construction detailed in Subsection~\ref{sec:cbc} for
any prime power $n$, and for each lattice point we solve the approximate
elliptic problem~\eqref{eq:weak} using one common finite element discretization in the
domain $D$. If $p=1$, then we assume in addition that~\eqref{eq:pone}
holds. Under the assumptions \textnormal{(A1)--(A6)}, we have
the combined error estimate
\begin{align*}
&\bigg|\int_UG(u(\cdot,\boldsymbol{y}))\,{\rm d}\bsy-\frac{1}{n} \sum_{k=0}^{n-1} G(u_h^s(\cdot,\boldsymbol{y}_k-\tfrac{\mathbf{1}}{\mathbf{2}}))\bigg|\\
&\leq C(\kappa(s,n)\|G\|_{H^{-1}}\|f\|_{H^{-1}}+h^{t+t'}\|G\|_{H^{-1+t'}}\|f\|_{H^{-1+t}}),
\end{align*}
where $0\leq t+t'\leq 2$, $h$ denotes the mesh size of the piecewise linear finite element mesh, $C>0$ is a constant independent of $s$, $h$, $f$, and $G$, and
$$
\kappa(s,n)=\begin{cases}
s^{-2/p+1}+n^{-1/p}&\text{if }p\in (0,1),\\
\big(\sum_{j=s+1}^\infty b_j\big)^2+n^{-1}&\text{if }p=1.
\end{cases}
$$
\end{theorem}

\section{Numerical experiments}\label{sec:numex}

We solve~\eqref{eq:weak} in the two-dimensional physical domain
$D=(0,1)^2$ with the source term $f(\bsx) = x_2$ and the periodic
diffusion coefficient~\eqref{eq:rfield}, denoted below by $a_{\rm
per}(\bsx,\bsy)$, where $\overline{a}(\bsx)= 2$ and
\begin{align}
\psi_{j}(\bsx)=c\,j^{-\beta}\sin(j\pi x_1)\sin(j\pi x_2)\quad\text{for}~c>0,~\beta>1~\text{and}~j\in\mathbb{N}.\label{eq:fluctuation}
\end{align}

For the numerical experiments, we truncate the parametric dimension to
$s=100$ and use a first order finite element solver to compute solutions
to~\eqref{eq:weak} 
numerically by using a regular finite element mesh of the square domain $D$ with the
one-dimensional mesh width $h=2^{-7}$. We use lattice rules generated
by the fast CBC algorithm detailed in Subsection~\ref{sec:cbc} with
$$n\in\{17,31,67,127,257,503,1\,009,2\,003,4\,001,8\,009,16\,007,32\,003,64\,007\}$$ nodes and choose $\sigma=\alpha=\beta\in\{2,4\}$. Moreover, all computations have been carried out using three different values for the scaling parameter $c\in\{1,0.5,0.1\}$ to allow us to vary the difficulty of the resulting integration problem. The reference solution was computed using a rank-1 lattice rule with $n=128\,021$
nodes.

In addition, we compare the convergence rates obtained in the periodic setting to the rates obtained using interlaced polynomial lattice rules generated for the
problem~\eqref{eq:weak}, equipped instead with the {\em affine}
diffusion coefficient
\begin{equation}
 %a(\bsx,\bsy)=
 a_{\rm aff}(\bsx,\bsy)=\overline{a}(\bsx)+\sum_{j\geq 1} y_j\,\psi_{j}(\bsx),\quad \bsx\in D,~\bsy\in U,\label{eq:toydiff2}
\end{equation}
which has the same mean field and covariance as the periodic field $a_{\rm per}$ when the fluctuations are chosen as in~\eqref{eq:fluctuation}. To generate interlaced polynomial lattice rules tailored for the affine diffusion coefficient, we used the QMC4PDE toolbox~\cite{qmc4pde,qmc4pde2} with the interlacing factors chosen to be equal to $\beta$ and $n=2^k$, $k\in\{4,\ldots,16\}$. In this case, the reference solution was computed using a corresponding interlaced polynomial lattice rule with $n=2^{17}$ nodes.

The quantity of interest in the first numerical experiment is the expectation $\mathbb{E}[G(u)]$ of the linear functional
$$
G(u)=\int_Du(\bsx,\bsy)\,{\rm d}\bsx,\quad \bsy\in U,
$$
\begin{figure}[!t]
\centering
\subfloat
{{\includegraphics[height=.429 \textwidth]{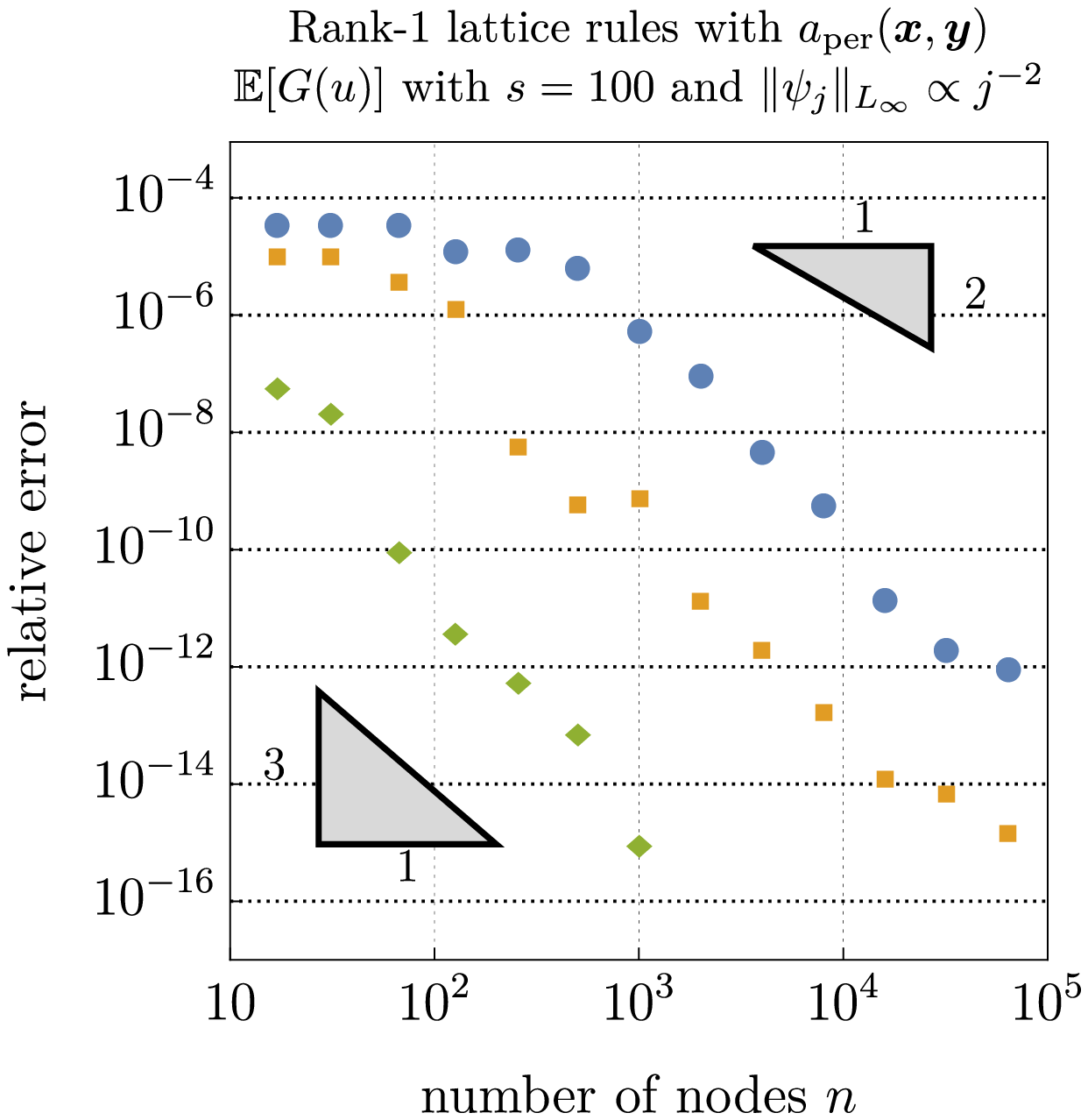}}}
\subfloat
{{\includegraphics[height=.429 \textwidth]{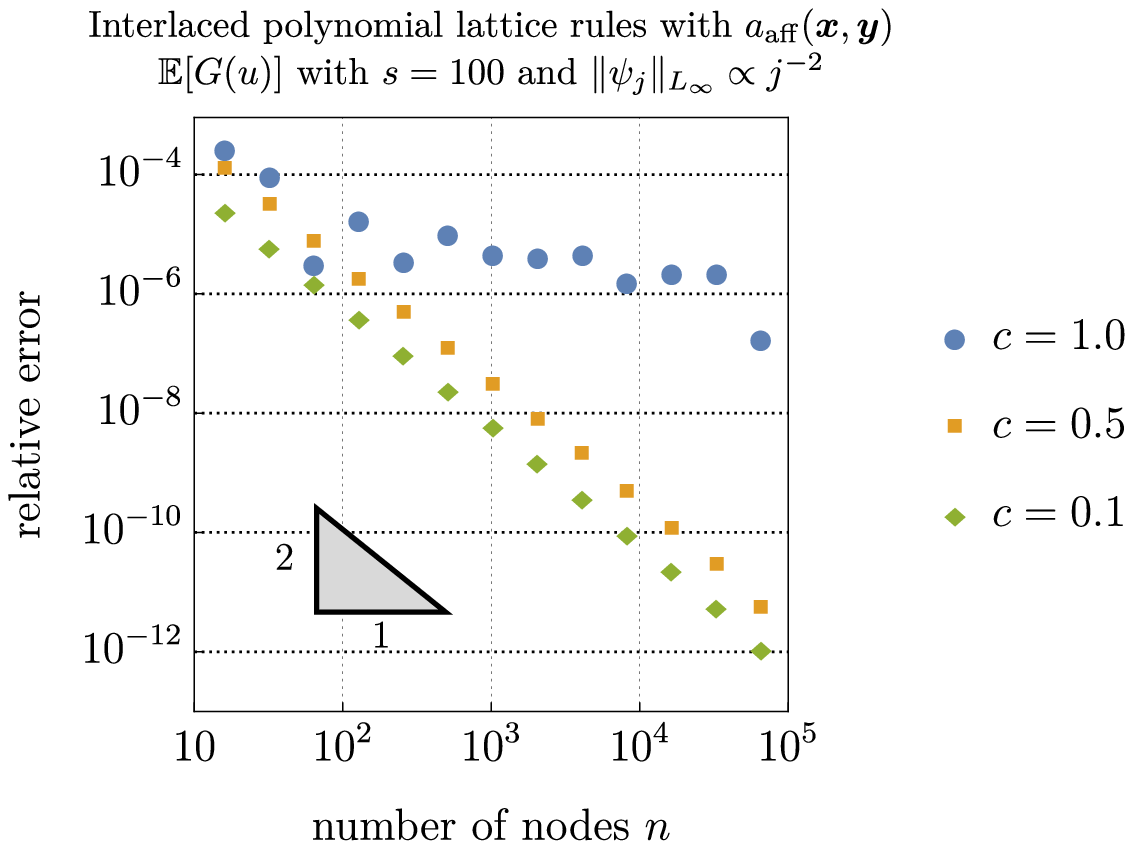}}}
\caption{Comparison of cubature errors in approximating $\mathbb{E}[G(u)]$ between rank-1 lattice rules in the periodic model~\eqref{eq:rfield} and interlaced polynomial lattice rules in the affine model~\eqref{eq:toydiff2} for $\beta=2$.}\label{fig:comparison1}
\end{figure}
\begin{figure}[!t]
\centering
\subfloat
{{\includegraphics[height=.429 \textwidth]{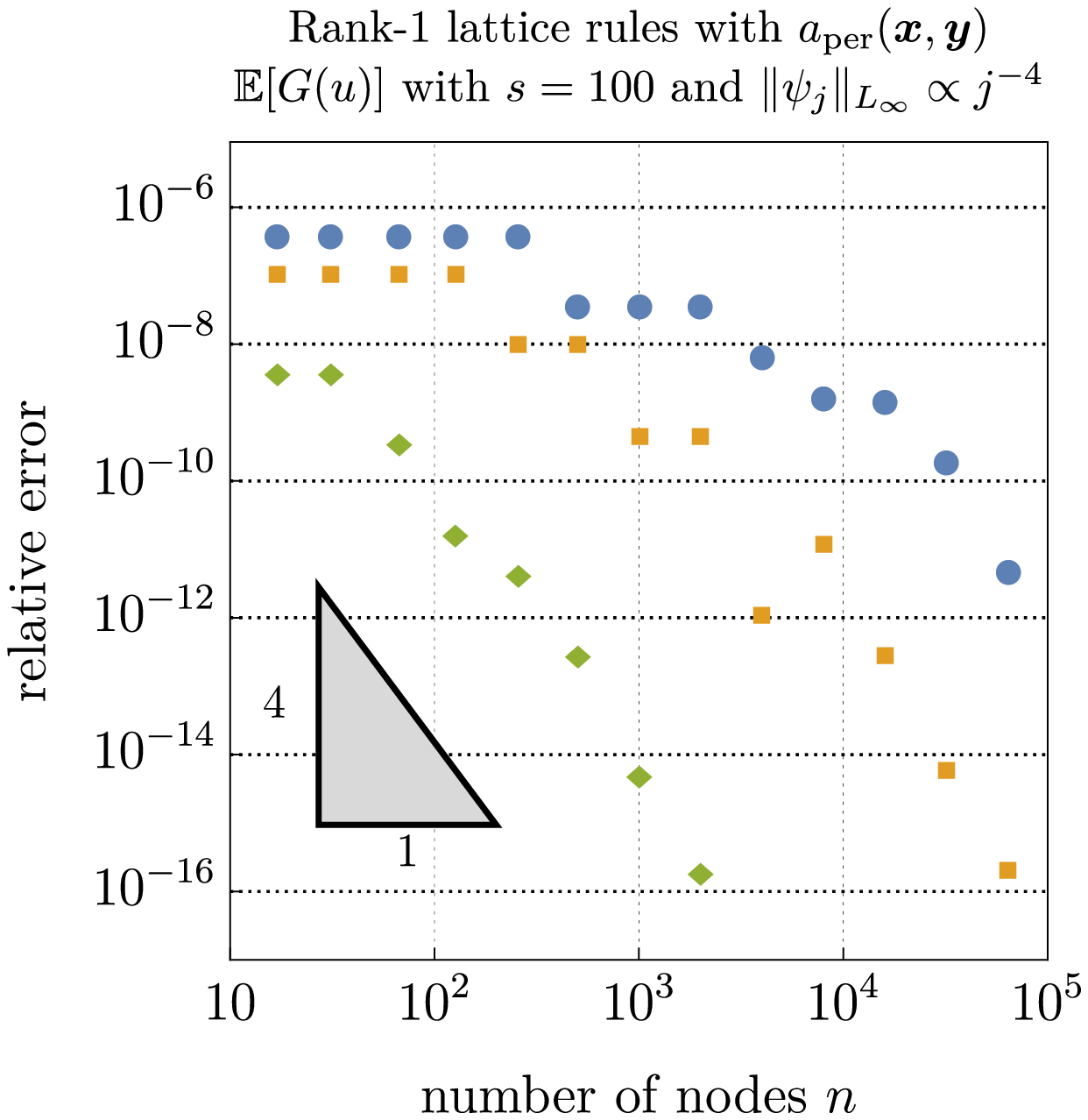}}}
\subfloat
{{\includegraphics[height=.429 \textwidth]{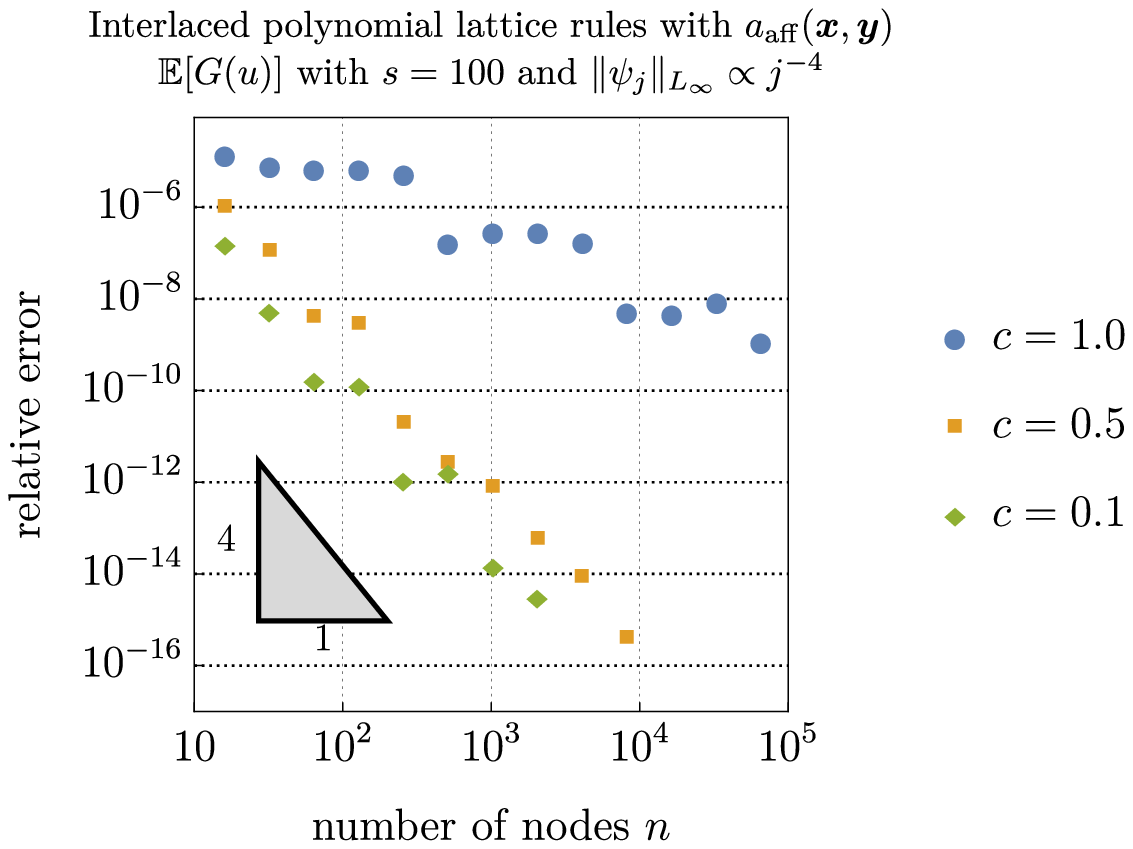}}}
\caption{Comparison of cubature errors in approximating $\mathbb{E}[G(u)]$ between rank-1 lattice rules in the periodic model~\eqref{eq:rfield} and interlaced polynomial lattice rules in the affine model~\eqref{eq:toydiff2} for $\beta=4$.}\label{fig:comparison1b}
\end{figure}
where the value of this integral can be calculated exactly when the integrand is a finite element solution. The results obtained using rank-1 lattice rules for the periodic
model~\eqref{eq:rfield} are displayed on the left-hand sides of
Figures~\ref{fig:comparison1} and~\ref{fig:comparison1b} for the decay rates $\beta=2$ and $\beta=4$, respectively. The corresponding results obtained using interlaced
polynomial lattices for the affine model~\eqref{eq:toydiff2} are  displayed on the
right-hand sides of Figures~\ref{fig:comparison1} and~\ref{fig:comparison1b}. The expected rates of convergence are $\mathcal{O}(n^{-2})$ and $\mathcal{O}(n^{-4})$, respectively. We observe that the solutions computed using the periodic diffusion coefficient $a_{\rm per}$ appear to converge at a rate at least as good as the expected rate. When the scaling parameter is set to $c=1$, it is notable that the rank-1 lattice rules used in conjunction with the periodic model appear to outperform the solution computed using interlaced polynomial lattice rules within the affine framework. It is apparent from Figure~\ref{fig:comparison1} that the observed rate of convergence is actually slightly better than the expected rate. This may be attributed to the fact that the dependence of the solution $u$ to the parametric variable $\bsy$ is analytic. For the solutions obtained using the affine model in Figure~\ref{fig:comparison1}, the interlacing factor $2$ actually acts as a bottleneck, capping the convergence rate at $\mathcal{O}(n^{-2})$. Similar numerical behavior for interlaced polynomial lattice rules has been previously reported, e.g., in~\cite{gantner2}.

\begin{figure}[!t]
\centering
\subfloat
{{\includegraphics[height=.429 \textwidth]{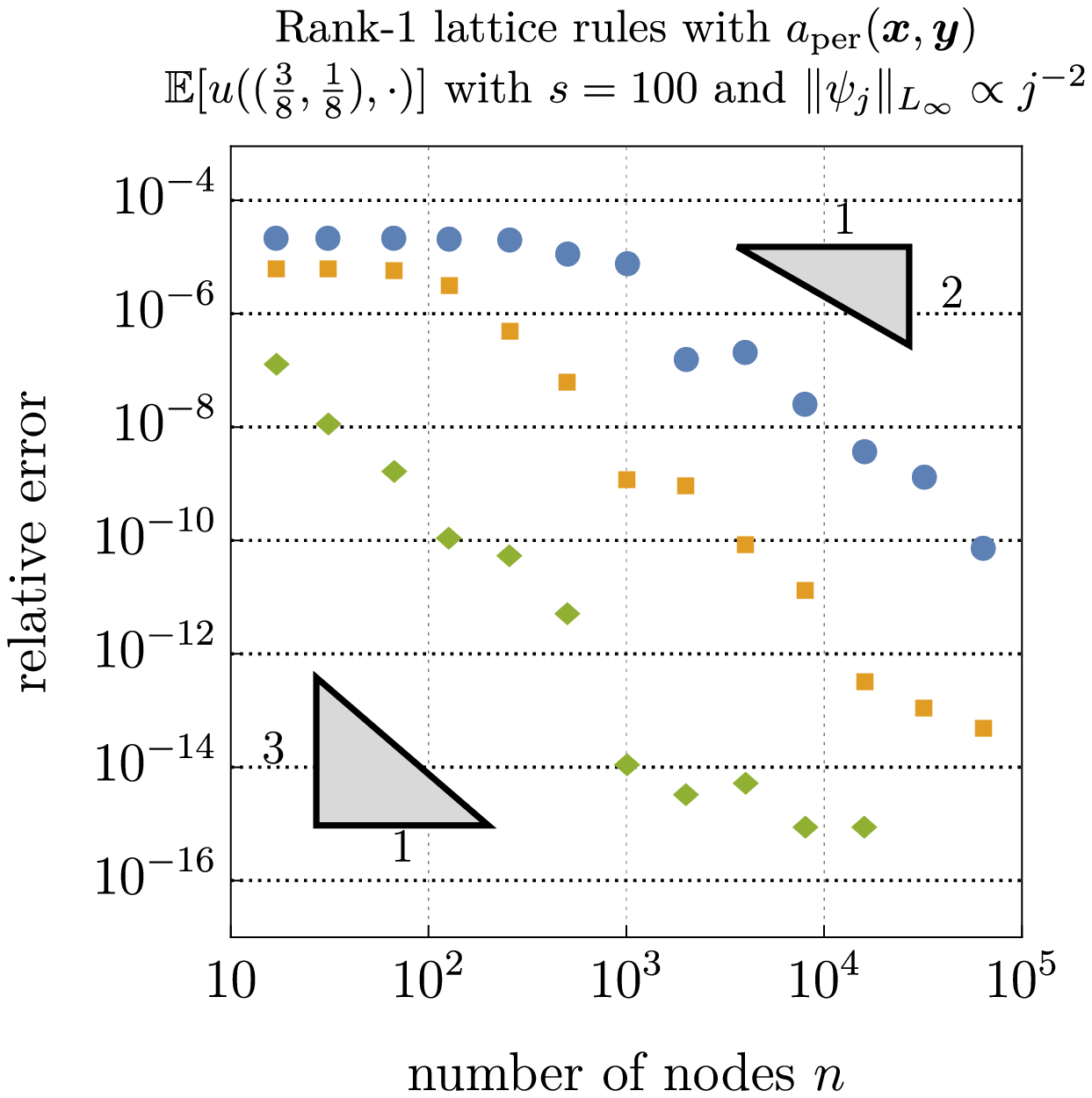}}}
\subfloat
{{\includegraphics[height=.429 \textwidth]{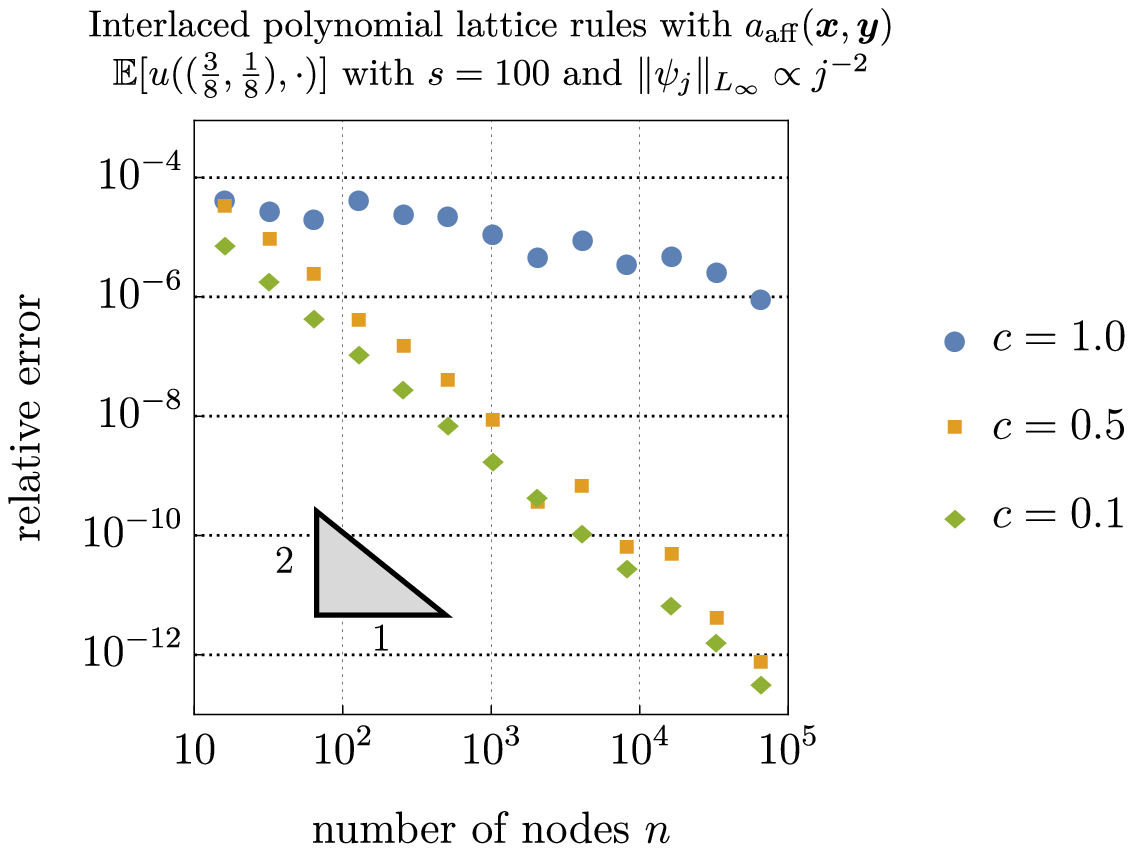}}}
\caption{Comparison of cubature errors in approximating $\mathbb{E}[u((\tfrac38,\tfrac18),\cdot)]$ between rank-1 lattice rules in the periodic model~\eqref{eq:rfield} and interlaced polynomial lattice rules in the affine model~\eqref{eq:toydiff2} for $\beta=2$.}\label{fig:comparison2}
\end{figure}
\begin{figure}[!t]
\centering
\subfloat
{{\includegraphics[height=.429 \textwidth]{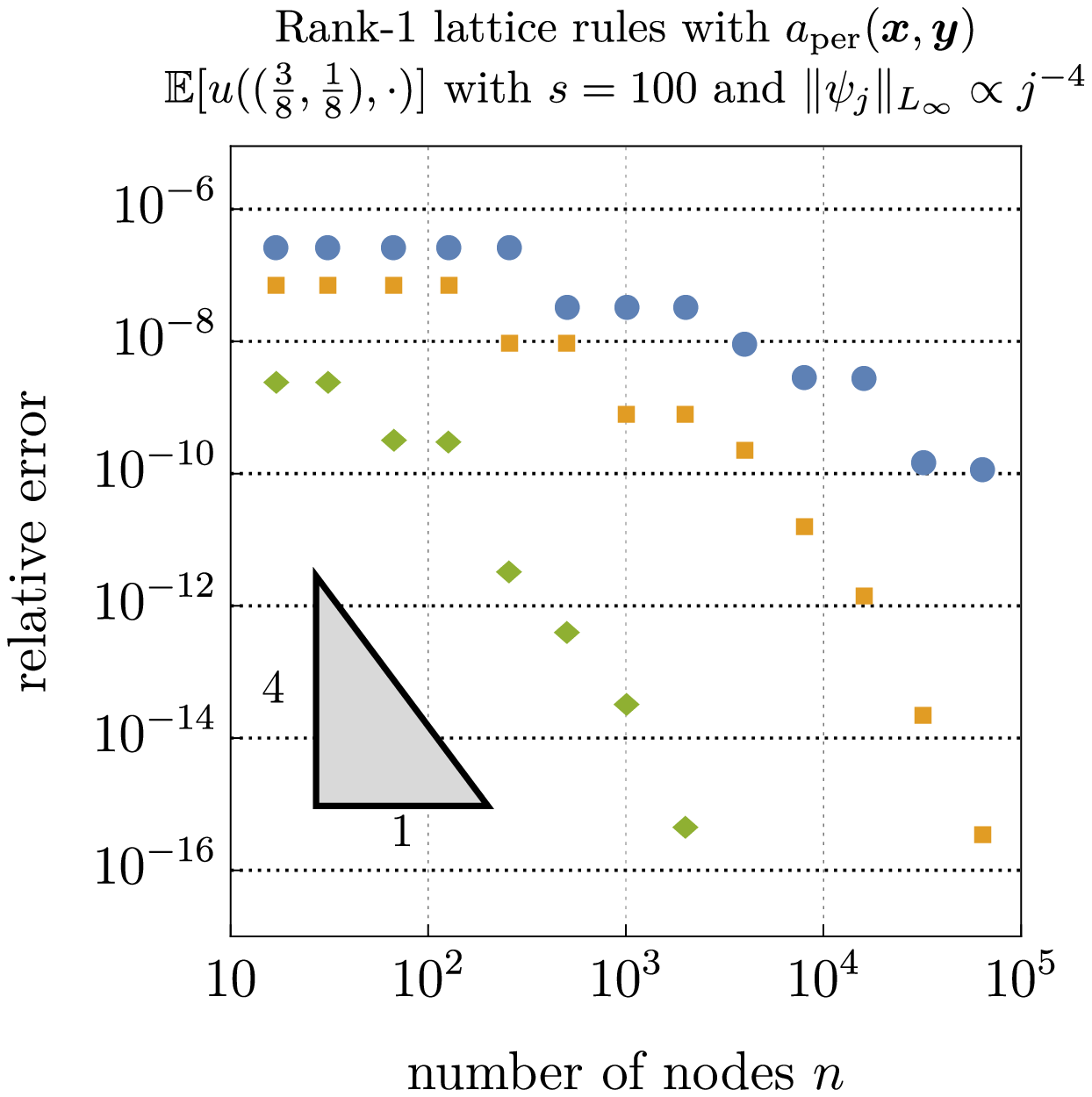}}}
\subfloat
{{\includegraphics[height=.429 \textwidth]{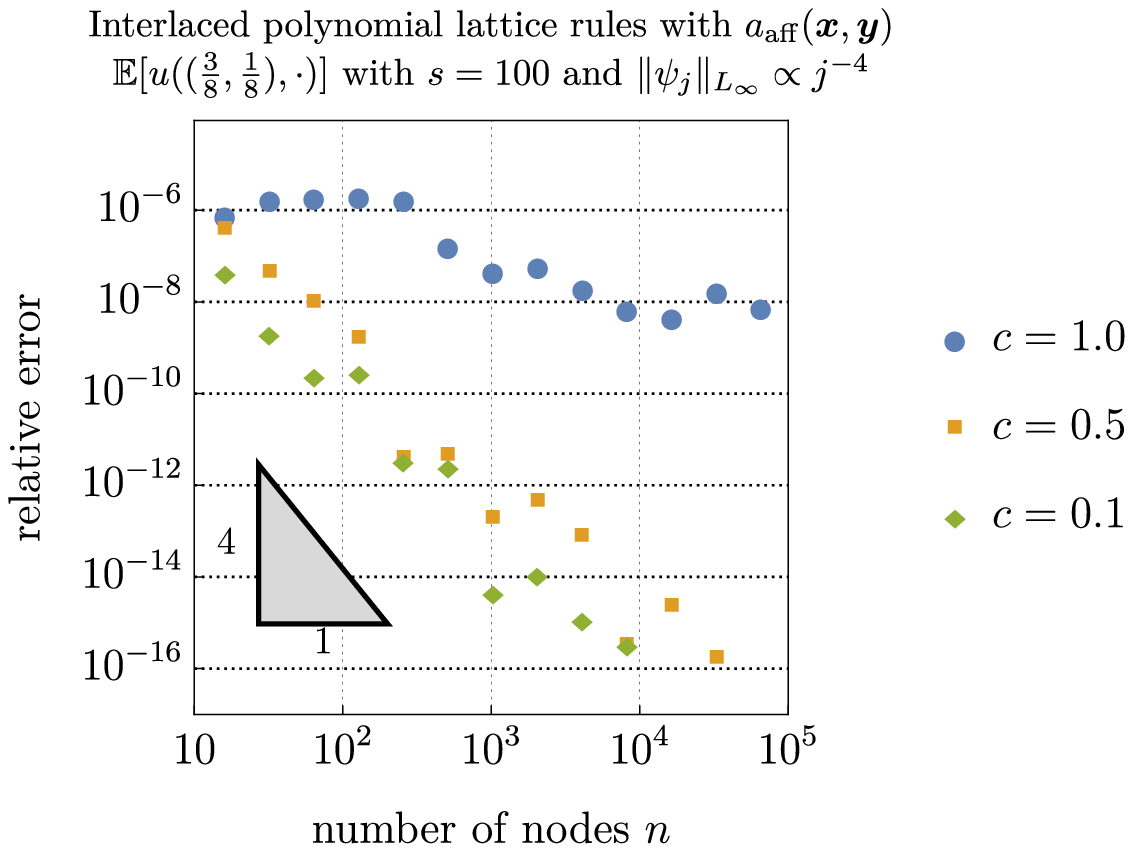}}}
\caption{Comparison of cubature errors in approximating $\mathbb{E}[u((\tfrac38,\tfrac18),\cdot)]$ between rank-1 lattice rules in the periodic model~\eqref{eq:rfield} and interlaced polynomial lattice rules in the affine model~\eqref{eq:toydiff2} for $\beta=4$.}\label{fig:comparison2b}
\end{figure}

In our second experiment, we consider the problem of approximating
$$
\mathbb{E}[u((\tfrac38,\tfrac18),\cdot)]=\int_Uu((\tfrac38,\tfrac18),\bsy)\,{\rm d}\bsy.
$$
The parameters and weights used for the construction of the rank-1 lattice
rules with the periodic diffusion coefficient~\eqref{eq:rfield} as
well as the interlaced polynomial lattice rules generated for the affine
diffusion coefficient~\eqref{eq:toydiff2} are exactly the same as in the
first numerical experiment. The results are displayed in Figures~\ref{fig:comparison2} and~\ref{fig:comparison2b}. We find that the general trend of the results matches that of the first numerical experiment, with the observed rates being at least as good as the expected rates with the scaling parameters $c\in\{0.5,0.1\}$, while the results obtained for the periodic model with $c=1$ and $\beta=4$ appear to remain in the preasymptotic regime.

{\em Remark.} Since the higher order moments of the input random
fields~{\eqref{eq:rfield}} and~\eqref{eq:toydiff2} are in general
different, so are the corresponding solutions to the respective
integration problems. Making a direct numerical comparison of the values
obtained in either setting is therefore not sensible.

\section*{Conclusions}
From a modeling point of view, there does not seem to be a reason to prefer an affine expansion of a  random field over a periodic expansion. Yet in the context of uncertainty quantification for PDEs with uncertain coefficients, we have seen that the model chosen for the random coefficient can make all the difference between obtaining essentially linear convergence with the affine model on the one hand, and on the other hand higher order convergence with the periodic model using rank-1 lattice cubature rules for the task of approximating the response statistics of the system. Higher order convergence can also be obtained with the affine model using, for example, interlaced polynomial lattice rules, but the overwhelming simplicity of constructing rank-1 lattice cubature rules makes the periodic framework a very enticing model for solving PDE problems equipped with uncertain coefficients. We have also presented numerical experiments that assess the QMC error derived in this work, in which the results are at least as good as those for a comparable affine model with interlaced polynomial lattice rules.
\section*{Acknowledgements}
We gratefully acknowledge the financial support from the Australian Research Council (DP180101356). We are also grateful to Fabio Nobile and Yoshihito Kazashi for collaboration on a related joint project that led to the inception of this paper. Frances Kuo thanks the participants from the Oberwolfach Workshop 1911 on Uncertainty Quantification for stimulating discussions about this work.

\bibliographystyle{plain}
\bibliography{UQ_periodic}
\end{document}